\documentclass{article}

\usepackage{PRIMEarxiv}

\usepackage[utf8]{inputenc} 
\usepackage[T1]{fontenc}    
\usepackage{hyperref}       
\usepackage{url}            
\usepackage{booktabs}       
\usepackage{amsfonts}       
\usepackage{nicefrac}       
\usepackage{microtype}      
\usepackage{lipsum}
\usepackage{fancyhdr}       
\usepackage{graphicx}       
\usepackage{subcaption}     
\usepackage{mathtools}
\usepackage{pgfplots}       
\pgfplotsset{compat=newest}

\DeclareMathOperator*{\argmin}{argmin}

\usepackage{tikz}
\usetikzlibrary{cd}

\usepackage[boxruled,linesnumbered,lined]{algorithm2e}
\usepackage[skins]{tcolorbox}
\newcommand{\Ke}[1]{\textcolor{blue}{Ke: #1}}

\newcommand{\Xcal}{\mathcal{X}}
\newcommand{\Ycal}{\mathcal{Y}}
\newcommand{\Scal}{\mathcal{S}}
\newcommand{\Ical}{\mathcal{I}}
\newcommand{\Acal}{\mathcal{A}}
\newcommand{\Fcal}{\mathcal{F}}
\newcommand{\Ocal}{\mathcal{O}}
\newcommand{\Lcal}{\mathcal{L}}
\newcommand{\Bcal}{\mathcal{B}}
\newcommand{\Mcal}{\mathcal{M}}
\newcommand{\Rcal}{\mathcal{R}}

\newcommand{\Rbb}{\mathbb{R}}
\newcommand{\Cbb}{\mathbb{C}}
\newcommand{\Ebb}{\mathbb{E}}
\newcommand{\Sbb}{\mathbb{S}}

\newcommand{\asf}{\mathsf{a}}
\newcommand{\vsf}{\mathsf{v}}
\newcommand{\bsf}{\mathsf{b}}
\newcommand{\usf}{\mathsf{u}}

\newcommand{\Psf}{\mathsf{P}}
\newcommand{\Qsf}{\mathsf{Q}}
\newcommand{\Vsf}{\mathsf{V}}
\newcommand{\Asf}{\mathsf{A}}
\newcommand{\Usf}{\mathsf{U}}
\newcommand{\Bsf}{\mathsf{B}}

\title{pseudo-differential integral autoencoder network for inverse PDE operators
}

\author{%
	Ke Chen\\
	Department of Mathematics \\
	University of Maryland College Park \\
	4176 Campus Dr, College Park, MD, 20742,  USA\\
	\texttt{kechen@umd.edu} \\
	\And 
	Jasen Lai\\
	Department of Statistics\\
	Purdue University\\
	150 N. University Street, West Lafayette, IN 47907\\
	\texttt{lai216@purdue.edu}
	\And
	Chunmei Wang \\
	Department of Mathematics \\
	University of Florida \\
	1400 Stadium Rd, Gainesville, FL, 32611, USA\\
	\texttt{chunmei.wang@ufl.edu}
	\\
}

\begin{document}
	\maketitle

	\begin{abstract}
		Partial differential equations (PDEs) play a foundational role in modeling physical phenomena. This study addresses the challenging task of determining variable coefficients within PDEs from measurement data. We introduce a novel neural network, "pseudo-differential IAEnet" (pd-IAEnet), which draws inspiration from pseudo-differential operators. pd-IAEnet achieves significantly enhanced computational speed and accuracy with fewer parameters compared to conventional models.
		Extensive benchmark evaluations are conducted across a range of inverse problems, including Electrical Impedance Tomography (EIT), optical tomography, and seismic imaging, consistently demonstrating pd-IAEnet's superior accuracy. Notably, pd-IAEnet exhibits robustness in the presence of measurement noise, a critical characteristic for real-world applications. An exceptional feature is its discretization invariance, enabling effective training on data from diverse discretization schemes while maintaining accuracy on different meshes.
		In summary, pd-IAEnet offers a potent and efficient solution for addressing inverse PDE problems, contributing to improved computational efficiency, robustness, and adaptability to a wide array of data sources.
		
	\end{abstract}

	
	
	\section{Introduction}
	
	Partial differential equations (PDEs) serve as foundational mathematical tools in various engineering and scientific computing domains. They play a crucial role in modeling and simulating a wide range of physical phenomena.
	The primary objective in tackling the "forward problem" is to determine the PDE solution when provided with information about the variable coefficients and initial/boundary conditions. This process is fundamental for understanding and predicting the behavior of physical systems.
	However, in many practical applications, the focus shifts towards the "inverse problem." This entails the challenging task of deducing the variable coefficients of the PDE based on multiple measurements of its solutions~\cite{isakov2006inverse}.  Solving inverse problems has long been a complex undertaking within the scientific computing community. Traditional approaches address nonlinear PDE inverse problems by formulating them as regression tasks over substantial datasets~\cite{chavent2010nonlinear}. These methods employ gradient descent techniques to seek local minimizers that represent the reconstructed coefficients. This iterative procedure demands numerous iterations, each involving thousands of forward and adjoint PDE solves, resulting in significant computational requirements. 
	One notable complication of inverse problems is their "ill-posed" nature. The candidate space for the coefficient function is exceedingly vast, and measurement errors can lead to substantial inaccuracies in the reconstruction. As such, it becomes imperative to incorporate prior knowledge and employ appropriate regularization techniques to confine the search space and enhance computational efficiency. Data-driven methods for PDE inverse problems have been an emergent field~\cite{arridge2019solving}.
	
	Neural networks have gained significant popularity as a versatile approach for approximating solutions to partial differential equations (PDEs). Their appeal lies in their ability to handle nonlinearity effectively and potentially mitigate the curse of dimensionality, a challenge that traditional methods often struggle with. Typical neural network methods include but are not limited to  Deep Ritz method\cite{yu2018deep}, Physics informed neural networks (PINNs) \cite{raissi2019physics}, DeepONet~\cite{DeepONet}, and other related methods\cite{chen2023friedrichs,brandstetter2022message}. Beyond learning PDE solutions, deep learning methods have been proposed for learning PDE operators~\cite{FNO,ong2022integral,kovachki2021neural,pathak2022fourcastnet,lu2022comprehensive}. PDE Operator learning methods involve discretizing the computational domain and representing functions as vectors on a discrete mesh. A deep neural network is then used to approximate the nonlinear mapping between these finite-dimensional spaces. Despite their success in various applications \cite{lin2021operator,cai2021deepm}, a limitation of these methods is that they often require retraining when solving problems on different meshes. To address these issues, methods with discretization-invariant properties have been introduced in recent research such as \cite{FNO,lu2022comprehensive,ong2022integral,ong2022integral}, offering flexible input and output discretization. Additionally, addressing the curse of dimensionality has been a critical challenge~\cite{stone1982optimal}. This curse implies that a substantial amount of data is needed to train a neural network effectively for learning an operator, particularly when dealing with high-dimensional input spaces. This issue has been investigated for finite dimensional operator learning \cite{bauer2019deep,nakada2020adaptive,schmidt2020nonparametric} and infinite dimensional  operator  learning \cite{liu2022deep,lanthaler2022error,chen2023deep}. 
	In the context of inverse PDE operator learning, unique numerical challenges arise. Unlike forward PDE operators, which typically map input functions to PDE solutions, inverse PDE operators map measurement data to reconstructed target parameter functions. This measurement data often includes pairs of boundary conditions and boundary measurements obtained from a series of experiments. These challenges become particularly severe because the input variable for the inverse PDE operator is the product of two boundary functions, leading to extremely high input dimensions. Consequently, it is crucial to design neural networks that can efficiently handle large input dimensions while maintaining a compact size. 
	
	In our research, we introduce a novel neural network architecture called "pseudo-differential IAEnet" (pd-IAEnet) designed specifically for learning the inverse PDE operator. Our work brings forth several key contributions:
	\begin{itemize}
		\item  \textbf{Inspired by Pseudo-Differential Operators:} We draw inspiration from the concept of pseudo-differential operators to create a unique neural network architecture tailored for solving inverse problems by approximating the inverse PDE operator. pd-IAEnet builds upon the foundation of IAEnet \cite{ong2022integral}, but it stands out due to its improved computational speed and enhanced accuracy. This improvement is achieved through a novel low-rank factorization neural network structure. Furthermore, we demonstrate that pd-IAEnet exhibits superior space complexity by achieving similar or even better accuracy with significantly fewer trainable parameters when compared to other baseline models.
		\item \textbf{Extensive Benchmark Evaluation:} We thoroughly assess the performance of pd-IAEnet across a diverse set of benchmark inverse problems. These include tasks such as solving the Calder\'on problem for Electrical Impedance Tomography (EIT), reconstructing scattering coefficients in optical tomography, and addressing inverse scattering problems in seismic imaging. Comparative evaluations against established baseline models such as FNO \cite{FNO}, DeepONet \cite{DeepONet}, ResNet, and IAEnet \cite{ong2022integral} consistently show that our proposed method achieves the highest accuracy across most scenarios.
		\item \textbf{Robustness to Measurement Noise:} We conduct a comprehensive examination of the robustness of pd-IAEnet to measurement noise in both training and testing datasets. Even when trained on noisy data, pd-IAEnet demonstrates a remarkable resilience to noise, maintaining a high level of accuracy when applied to noisy testing data for inverse scattering problems. This robustness is crucial in real-world scenarios where measurement data often contains noise. For highly ill-posed inverse problems like optical tomography and the Calder\'on problem, pd-IAEnet achieves the best accuracy among all base models.
		\item \textbf{Discretization Invariance:} Another notable feature of pd-IAEnet is its discretization invariance property. This means that our method can be trained on measurement data obtained from various discretization schemes, while still preserving the same level of accuracy when applied to measurement data on different meshes. This flexibility is advantageous in practical applications where data may come from different sources.
	\end{itemize}
	
	This paper is organized as follows. We shall first introduce the mathematical formulation of PDE inverse problems and list a few benchmark inverse problems in Section \ref{sec:ip}. We shall introduce the pd-IAE net and discuss motivations, architectures and miscellaneous features in Section \ref{sec:pdiae}. Finally in Section \ref{sec:numerics}, we shall compare the proposed method with other baseline models on several benchmark inverse problems.  
	\section{Inverse problems}\label{sec:ip}
	\subsection{General mathematical formulation}
	Consider a general boundary value problem that seeks $u$ such that
	\begin{equation}\label{eqn:forward}
		\begin{cases}
			\Lcal_a u = h, \quad \textbf{in  } \Omega \subset \Rbb^d, \\
			\Bcal u = f, \quad \textbf{on  } \partial \Omega \,,
		\end{cases}
	\end{equation}
	where  $a:\Omega \rightarrow \Rbb$ is the unknown target function, $\Lcal_a$ denotes a differential operator defined within the domain $\Omega$, $h:\Omega \rightarrow \Rbb$ is a  known source function, $\Bcal$ is an operator defined on the boundary $\partial\Omega$, and $f: \partial \Omega \rightarrow \Rbb$ is the boundary condition. 
	
	\textbf{Inverse Problem and Cauchy Data:} In practical applications, multiple experiments are conducted with various boundary sources and receivers. These experiments are modeled using \eqref{eqn:forward}, where different boundary conditions $f$ lead to corresponding boundary measurements $g$, defined as $g \coloneqq \Mcal u$. The primary objective of the inverse problem is to reconstruct the unknown function $a$ based on information obtained from the "forward operator" $\Lambda_a: f\mapsto g$. This is equivalent to working with "Cauchy data":
	\begin{equation*}
		\Scal_a \coloneqq \{(f,g) \ |\ f\in C^\infty(\partial \Omega)\,, g = \Lambda_a f\} \,.
	\end{equation*}
	In practical scenarios, only a finite number of experiments can be conducted, resulting in a finite dataset:
	\begin{equation}\label{eqn:finitedata}
		\Scal_a^n \coloneqq \{(f_i,g_i) \ |\ f\in C^\infty(\partial \Omega)\,, g = \Lambda_a f\,, i = 1,\ldots,n\} \,.
	\end{equation}
	Here, the functions $f_i$ are predetermined and depend on the experimental design for configuring sources, while the measurements $g_i$ depend on the setup of receivers.

	\textbf{Traditional Computational Approach:} Typically, the inverse problem is formulated as the following regression problem~\cite{engl1996regularization}
	\begin{equation}\label{form:optimization}
		\hat{a} \in \argmin_a \frac{1}{n} \sum_{i=1}^n\| g_i - \Lambda_a f_i\|^2 + \Rcal(a) \,,
	\end{equation}
	where a regularization term $\Rcal(a)$ is introduced to incorporate prior knowledge about $a$. The optimization formulation \eqref{form:optimization}, assuming solvability, generally yields a unique solution when the number $n$ of measurement data is sufficiently large. However, solving this problem computationally poses challenges due to the nonlinearity of $\Lambda_a$. This necessitates iterative solvers, where each iteration involves solving the forward PDE \eqref{eqn:forward} and the adjoint PDE for all data pairs in $\Scal_a^n$. This computational approach becomes inefficient when $n$ is large, and selecting an appropriate regularization function $\Rcal$ can be complex.
	
	\textbf{Proposed Deep Learning Approach:} To address these computational challenges, we propose a deep learning approach. Specifically, we introduce an \textit{inverse operator} $\psi: S_a^n \mapsto a$ that maps finite data sets to the target parameter function. Our objective is to utilize a neural network $\phi_\theta$ to parameterize $\psi: S_a^n \mapsto a$. This neural network, characterized by its parameters $\theta$, is trained on multiple data pairs $(a_j,S_{a_j}^n)$ by optimizing these parameters. Once trained with an optimal parameter $\theta^\ast$, the neural network $\phi_{\theta^\ast}$ can be employed to estimate any target function $a$ by evaluating it with the corresponding measurement data $\Scal_a^n$. Further details about this neural network approach will be discussed in Section \ref{sec:pdiae}.
	
	
	\subsection{Benchmark inverse problems} In this section,  we shall  present three benchmark inverse problems.

	\subsubsection{Inverse scattering}
	The inverse scattering problem~\cite{cakoni2005qualitative} is typically described using the Helmholtz equation:
	\begin{equation}\label{eqn:helmholtz}
		\left(- \Delta - \frac{w^2}{c(x)^2}\right) u = 0 \,, \quad x \in \Rbb^d \,,
	\end{equation}
	where $w \in \Rbb_+$ represents the frequency and $c(x)>0$ is the unknown inhomogeneous wave speed. 
	Assuming an incident field  given by a plane wave $u^i(x) = e^{iw s\cdot x}$ with a source direction $s\in \Sbb^1$, the total field can be written as $u(x) = e^{iw s\cdot x} + u^s(x)$ where $u^s(x)$ is the scatter field. The scatter field satisfies the Sommerfeld radiation condition:
	\[
	\lim_{r\rightarrow \infty} r\left( \frac{\partial u}{\partial r} - i w u \right) = 0 \,.
	\]
	Here $r = \| x\|$ and the limit is assumed to hold uniformly in all directions $\frac{x}{\|x\|}$. This condition implies that the scatter field is an outgoing spherical wave:
	\[
	u^s(x) = \frac{e^{iw\|x\|}}{\|x\|} \left(u^s_\infty(\frac{x}{\|x\|}) + o(\frac{1}{\|x\|})\right) \,,\quad \|x\| \rightarrow \infty \,.
	\]
	Here $u^s_\infty$ is a function defined on the unit sphere, known as the far field pattern of $u^s$. Measurements of the far field pattern can be collected on multiple receiver directions $r\in\Sbb^1$. More details of this model can be found in \cite{colton1998inverse}.
	
	In many applications, it is assumed that a background wave speed $c_0(x)$ is known and identical to $c(x)$ except on a bounded domain $\Omega$ within the unit ball. The scatter function $\eta(x)$ is then defined as:
	\begin{equation}\label{eqn:scatter}
		\eta(x) = \frac{w^2}{c(x)^2} - \frac{w^2}{c_0(x)^2} \,, \quad x\in \Omega \,.
	\end{equation}
	As a result, the scatter function $\eta$ has bounded support $\Omega$.
	
	The objective of the inverse scattering problem is to reconstruct the scatter function $\eta(x)$ from measurements of the far-field pattern  $d(r,s) = u^s_\infty(r)$  collected at multiple source and receiver locations $(s,r)\in \Sbb^1\times \Sbb^1$.  The well-posedness of inverse scattering problems has been extensively studied and reviewed in \cite{isakov1992stability,bao2015inverse}.

	\subsubsection{Optical tomography}
	Optical tomography (OT) plays a crucial role in reconstructing optical medium properties based on measurements of light transmitted and scattered through the medium, with significant applications in biomedical imaging for tissues like the brain and breast~\cite{arridge1999optical}. OT is mathematically modeled using the radiative transfer equation (RTE):
	\begin{equation}\label{eqn:rte}
		\begin{cases}
			v \cdot \nabla \rho(x,v)  = \sigma_s(x) \left( \int_{\Sbb^{d-1}} u(x,v')dv' - u(x,v) \right) \,, \quad (x,v) \in \Omega \times \Sbb^{d-1}\,, \\
			\rho(x,v) = f(x,v)\,,\quad (x,v)\in \Gamma_- \,,\\
		\end{cases}
	\end{equation}
	where $\sigma_s(x) >0$ represents the scattering coefficient  which characterizes the optical properties of the medium, and $\Gamma_-=\{ (x,v)\in \partial\Omega \times \Sbb^{d-1} \ |\ n_x \cdot v < 0\}$  defines the incoming boundary, where the incoming boundary condition $f(x,v)$ represents sources of photons injected into the domain. The outgoing photon intensity is denoted as $g\coloneqq \rho(x,v)|_{\Gamma_+}$ and is measured on the outgoing boundary $\Gamma_+=\{ (x,v)\in \partial\Omega \times \Sbb^{d-1} \ |\ n_x \cdot v > 0\}$. The primary objective of OT is to reconstruct the scattering coefficient $\sigma_s$ from multiple pairs of incoming and outgoing photon intensities, represented as $S_{\sigma_s}^n = \{(f_i,g_i)\ |\ i=1,\ldots,n\}$. OT is an ill-posed inverse problem that has been extensively studied, with comprehensive reviews available in references such as \cite{arridge2009optical,bal2009inverse}. 
	
	\subsubsection{Electrical impedance tomography}
	Electrical impedance tomography (EIT)~\cite{borcea2002electrical,cheney1999electrical}, often referred to as the Calderón problem~\cite{calderon1980inverse}, plays a significant role in noninvasive medical imaging, particularly in the early diagnosis of breast cancer. EIT aims to determine the electrical conductivity distribution within a medium by analyzing voltage and current measurements acquired along its boundary. Mathematically, EIT is formulated as an elliptic equation with Dirichlet boundary conditions:
	\begin{equation}\label{eqn:calderon}
		\begin{cases}
			- \text{div}\left(e^{a(x)} \nabla u(x)\right) = 0 \,, \quad x\in\Omega \,, \\
			u(x)= f(x) \,,\quad x \in \partial \Omega \,,
		\end{cases}
	\end{equation}
	where the term $e^{a(x)}>0$ represents the unknown medium conductivity, and $f$ denotes the boundary voltage. Current measurements, crucial for EIT, are defined as the Neumann derivatives of the solution:  $g \coloneqq e^{a} \frac{\partial u}{\partial n}|_{\partial\Omega}$ along the boundary. The main objective of the Calder\'on problem is to reconstruct the function $a(x)$, given the Dirichlet-to-Neumann data pairs $S_a^n = \{(f_i,g_i)\ |\ i = 1,\ldots, n\}$. The uniqueness of reconstructing the medium's conductivity has been established \cite{sylvester1987global} and the Calder\'on problem has been shown to be extremely ill-posed, with a logarithmic estimate in the stability result \cite{alessandrini1988stable}.
	
	\section{Pseudo-differential autoencoders algorithm}\label{sec:pdiae}
	
	The primary objective of the pseudo-differential integral autoencoders network (pd-IAEnet) is to learn an inverse PDE operator  $\Phi:\Xcal \rightarrow \Ycal$ using a neural network $\Phi_\text{NN}(\cdot;\theta_\text{NN}):\Xcal \rightarrow \Ycal$, where $\Xcal$  represents the measurement data space, $\Ycal$ is the space of target parameters, and $\theta_\text{NN}$ encompasses all trainable parameters in the neural network $\Phi_\text{NN}$.
	
	In the discrete setting, the finite measurement data \eqref{eqn:finitedata} is typically represented as:
	\[
	S^n_a = \left[S_{ij}\right] \in \Rbb^{n \times n}\,.
	\]
	Here $n$ is the number of grid points on the boundary, and each entry $S_{ij}$ corresponds to the measurement data collected at receiver $j$ for the solution generated by the $i$-th source on the boundary. In this discrete representation, the target medium function $a$ can be expressed as a vector $a = [a_i]\in \mathbb{R}^{M}$ that tabulates the function values $a_i = a(x_i)$, where $x_i$ represents the grid points, and $M$ is the total number of grid points in the domain. Consequently, the objective is to train a neural network to learn a nonlinear mapping from $\mathbb{R}^{N}$ to $\mathbb{R}^M$, where $N=n^2$. However, different experiments may employ distinct discretizations for the computational domain, leading to varying numbers of grid points $n$ on the boundary and $M$ in the domain. To address this challenge, the neural network must exhibit \textit{discretization invariance}, implying that its input and output should remain consistent regardless of the discretization.
	
	Operator learning is a burgeoning field in scientific machine learning that leverages neural networks to approximate nonlinear operators. The simplest neural network architecture is the fully-connected neural network (FNN), which consists of layers comprising dense affine transformations and pointwise nonlinear activation functions. A $N$ layer (or $N-1$ hidden layer) FNN $\Phi_\text{FNN}(x;\theta_\text{FNN}): \Rbb^{n_0} \rightarrow \Rbb^{n_N}$ is structured as follows:
	\begin{equation*}
		\Phi_\text{FNN}(x;\theta_\text{FNN}) = f_N \circ f_{N-1} \circ \cdots \circ f_1(x)\,,
	\end{equation*}
	where $f_i(x) = \sigma(W_i x + b_i)$ represents a composition of a pointwise activation function $\sigma$ and a linear transformation with weight matrix $W_i\in \Rbb^{n_i\times n_{i-1}}$ and bias vector $b_i\in \Rbb^{n_i}$, and $\theta_\text{FNN}$ encompasses all trainable parameters, including weight matrices $W_i$ and bias vectors $b_i$.  
	
	However, FNNs cannot be directly trained on datasets with varying data formats because both the input dimension $n_0$ and the output dimension $n_N$ are fixed. Additionally, FNNs cannot be applied to testing data with different discretizations, leading to incompatible input/output sizes. Furthermore, while neural networks like FNO\cite{FNO} and its variants can handle data of various discretizations, their accuracy may degrade when tested on data with discretizations different from the training data.
	
	To overcome these challenges, a novel integral autoencoder network (IAE-net) was introduced in \cite{ong2022integral} to achieve discretization invariance. However, the computation of IAE-net can be slow due to the neural network-parametrized integral kernel. Drawing inspiration from the pseudo-differential operator \eqref{eqn:def_pdo}, the pd-IAEnet employs a separation of variables in the integral kernel, resulting in faster evaluation and improved accuracy.
	
	\subsection{Pseudo-differential autoencoder networks (pd-IAEnet)}
	The pd-IAEnet follows a computational flow similar to the IAEnet and applies a sequence of discretization-invariant pd-IAE blocks $\Ical_1,\ldots,\Ical_L$ to maximize its approximation capabilities.
	
	Given an input data $\Bar{f}$ with a specific discretization $S_\Xcal$, the pd-IAEnet produces an output data $\Bar{g}$ with a given discretization $\Scal_\Ycal$ through the following computational flow:
	\begin{equation}\label{eqn:basic_IAE}
		\Bar{f} \xrightarrow{F} a_0 \xrightarrow{\Ical_1} a_1 \xrightarrow{\Ical_2}\cdots \xrightarrow{\Ical_L} a_L \xrightarrow{G} \Bar{g}
	\end{equation}
	where $F(\cdot;\theta_F)$  represents a pre-processing neural network function responsible for transforming the input data into a higher-dimensional space to enhance its features, and $G(\cdot;\theta_G)$  is a post-processing function serving a similar purpose. The data processing operators $F$ and $G$ utilized here are the same as those used in \cite{FNO}.

	\subsubsection{pd-IAE blocks and pd-encoders}
	The pd-IAE block is designed to map a function $a$ on a domain $\Omega_a$ to another function $b$ on domain $\Omega_b$, with both domains assumed to be $[0,1]^d$ for simplicity. The functions $a$ and $b$ are discretized on the same grid points $S = {x_i}_{i=1}^s$, which are allowed to vary.
	
	Each pd-IAE block consists of three components: an encoder function $a\mapsto v$, a FNN function $v\mapsto u$, and a decoder function $u\mapsto b$. In other words, it involves the following compositions:
	\begin{equation}\label{eqn:IAE-block}
		\Ical: \quad a \xrightarrow{\text{encoder}} v \xrightarrow{\text{FNN}} u \xrightarrow{\text{decoder}} b \,.
	\end{equation}
	The intermediate functions $v$ and $u$ are defined on $\Omega_z = [0,1]^d$ with a fixed number $m$ of grid points $S_z = \{z_j\}_{j=1}^m \subset \Omega_z$. The FNN within a pd-IAE block has fixed input and output dimensions of $m$, whereas the input and output dimensions $s$ for a pd-IAE block may vary for different training and testing data pairs.
	
	The encoder function in the original IAEnet uses a nonlinear integral transform with an NN-parametrized kernel $\phi_1(\cdot;\theta_{\phi_1})$. The encoder's output $v$ is computed via this integral transform
	\begin{equation}\label{eqn:iae-en}
		v(z) = \int_{\Omega_a} \phi_1\left(a(x),x,z;\theta_{\phi_1} \right) a(x) dx\,,\quad z \in \Omega_z \,.
	\end{equation}
	Analogously, the decoder function is also an integral transform 
	\begin{equation}\label{eqn:iae-de}
		b(x) = \int_{\Omega_z} \phi_2\left(u(z),x,z;\theta_{\phi_2} \right) u(z) dz\,,\quad x \in \Omega_b \,,
	\end{equation}
	with an NN-paramtrized kernel $\phi_2(u(z),x,z;\theta_{\phi_2})$. The nonlinear integral transform structure of the encoder and decoder is the key to discretization invariance of the IAE-net.
	These integral transforms depend on input grid points $x$, output grid points $z$, and the function values $a(x)$ (or $u(z)$). However, this dependence on input and output locations, as well as function values, can lead to computationally expensive encoders and decoders, and learning the kernels $\phi_1$ and $\phi_2$ can be challenging.

	The pd-IAE net simplifies the encoder function by introducing pseudo-differential integral autoencoders (pd-encoders) that resemble the formulation of a pseudo-differential operator $P: C_c^\infty(\Omega)\rightarrow C^\infty(\Omega)$. This operator is defined as:
	\begin{equation}\label{eqn:def_pdo}
		Pu(x) = \Fcal^{-1}(\Acal(x,\cdot) \hat{u}(\cdot))= \frac{1}{(2\pi)^n} \int e^{i x\cdot \xi} \Acal(x,\xi) \Hat{u}(\xi) d\xi \,,
	\end{equation}
	where $\Hat{u}$ is the Fourier transform of a $C_c^\infty(\Omega)$ function $u$, and $\Acal(x,\xi)$ is a smooth function that known as the symbol of $P$. These pseudo-differential operators generalize classic differential operators. For instance, if $\Acal(x,\xi)$ is a polynomial $p(\xi)$ in $\xi$, then the operator $P$  corresponds to the classical differential operator $p(-i\partial_x)$.
	
	Note that due to the multi-channel structure \eqref{eqn:multi-channel}, the input function $a(\xi)$ has already been transformed into the frequency domain through the Fourier transform. This simplifies the encoder function in the pd-IAEnet, which now takes the form:
	\begin{equation}\label{eqn:pd_encoder}
		v(z) =\Fcal^{-1}\left(\Acal(z,\cdot)a(\cdot)\right)= \int_{\Omega_a} e^{i z\cdot \xi }A(z,\xi) a(\xi)  d\xi \,.
	\end{equation}
	This simplified encoder function involves only one Fourier transform and does not depend on input grid points. To approximate the kernel function $A(z,\xi)$, we use a low-rank factorization:
	\begin{equation}\label{eqn:low_rank}
		A(z,\xi) \approx \sum_{i=k}^K p_k(\xi)q_k(z) \,.
	\end{equation}
	Combining \eqref{eqn:pd_encoder} and \eqref{eqn:low_rank}, the encoder function can be further simplified:
	\begin{equation}\label{eqn:low_rank_encoder}
		v(z) \approx \sum_{k=1}^K q_k(z)\int_{\Omega_a} e^{i z\cdot \xi }  a(\xi) p_k(\xi) d\xi \,.
	\end{equation}
	Inspired by this simplified low-rank encoder function, we propose to replace the functions $p_k$ and $q_k$ with neural network parametrized functions:
	\begin{equation}\label{eqn:low_rank_encoder2}
		v(z) \approx \sum_{k=1}^K q_k(z;\theta_{q_k})\int_{\Omega_a} e^{i z\cdot \xi }  a(\xi) p_k(\xi;\theta_{p_k}) d\xi 
		= \sum_{k=1}^K q_k(z;\theta_{q_k}) \Fcal^{-1} \left(p_k(\cdot;\theta_{p_k}) a(\cdot)\right)\,,
	\end{equation}
	where $p_k(\xi;\theta_{p_k})$ and $q_k(z;\theta_{q_k})$ are complex-valued functions defined over the domain  $\Cbb^1$. 
	Since the input function $a(\xi)$ is already in the frequency domain, we can safely truncate high-frequency modes in $a(\xi)$, ensuring computational efficiency while maintaining accuracy. This truncation is crucial for achieving discretization invariance in the discrete setting.
	
	In particular, the original input function $a(x)\Rbb^s$ consists of function values at arbitrary grid points $S = \{x_i\}_{i=1}^s$ and is mapped to a vector $a(\xi)$ with fixed size due high frequency truncation.
	This ensures that the intermediate FNN within an IAE-block has a fixed input size.
	Once all the bases functions $p_k$ and $q_k$  have been trained, they form two data matrices $\Psf$ and $\Qsf$. Similar to the pseudo-differential encoder function \eqref{eqn:low_rank_encoder2}, the pseudo-differential decoder function $u\mapsto b$ has the following structure
	\begin{equation*}
		b(x) = \sum_{k=1}^K \Tilde{p}_k(x;\theta_{\Tilde{p}_k}) \Fcal\left( \Tilde{q_k}(\cdot;\theta_{\Tilde{q}_k}) u(\cdot)\right) \,,
	\end{equation*}
	where $\Tilde{p}_k$  and $\Tilde{q}_k$ are parametrized NNs. 
	The Fourier transform and inverse Fourier transform operations can be efficiently implemented using fast Fourier transform (\texttt{fft}) and fast inverse Fourier transform (\texttt{ifft}) techniques. 
	
	To summarize, the computation of the encoder and decoder functions can be implemented using the procedures outlined in Algorithm \ref{alg:encoder} and Algorithm \ref{alg:decoder}, respectively.
	
	\vspace{3 pt} 
	\strut\vspace*{-\baselineskip}\newline
	\begin{minipage}{0.48\textwidth}
		\begin{algorithm}[H]
			\caption{pd-encoder function}\label{alg:encoder}
			\KwIn{$\asf \in \Rbb^{s\times 1}\,, m\leq s \in \mathbb{N}$}
			\KwData{data matrices $\Psf \in \Rbb^{m\times K}$ and $\Qsf \in \Rbb^{m\times K}$} 
			\KwOut{$\vsf \in \Rbb^{m\times 1}$}
			Truncate high frequency modes: $\asf = \asf(1:m)$\;
			Compute $\Asf =\asf \mathsf{1}_K^\top \in \Rbb^{m\times K}$\;
			Compute pointwise product $\Psf = \Psf \odot \Asf \in \Rbb^{m\times K}$\;
			\For{$i = 1:K$}{
				$\Psf(:,i) = \texttt{ifft}(\Psf(:,i))$\;
			}
			Compute pointwise product $\Vsf = \Qsf \odot \Psf \in \Rbb^{m\times K}$\;
			\For{$i = 1:m$}{
				$\vsf(i) = \texttt{sum}(\Vsf(i,:)) $ \;
			}
		\end{algorithm}
	\end{minipage}
	\hfill
	\vspace*{3 pt} 
	\begin{minipage}{0.48\textwidth}
		\begin{algorithm}[H]
			\caption{pd-decoder function}\label{alg:decoder}
			\KwIn{$\usf \in \Rbb^{m\times 1}\,, s \geq m \in \mathbb{N}$}
			\KwData{data matrices $\Tilde{\Psf} \in \Rbb^{s\times K}$ and $\Tilde{\Qsf} \in \Rbb^{s\times K}$} 
			\KwOut{$\bsf \in \Rbb^{s\times 1}$}
			Compute $\Usf =\Hat{\usf} \mathsf{1}_K^\top \in \Rbb^{m\times K}$\;
			Pad zeros $\Usf = \begin{bmatrix}
				\Usf \\
				\mathsf{0}
			\end{bmatrix} \in \Rbb^{s\times K}$\;
			Compute pointwise product $\Tilde{\Qsf}= \Tilde{\Qsf} \odot \Usf \in \Rbb^{s\times K}$\;
			\For{$i = 1:K$}{
				$\Tilde{\Qsf}(:,i) = \texttt{fft}(\Tilde{\Qsf}(:,i))$\;
			}
			Compute pointwise product $\Bsf =\Tilde{ \Qsf} \odot \Tilde{\Psf} \in \Rbb^{s\times K}$\;
			\For{$i = 1:s$}{
				$\bsf(i) = \texttt{sum}(\Bsf(i,:)) $ \;
			}
		\end{algorithm}
	\end{minipage}
	\\
	Note that the use of low-rank factorization in the integral kernel \eqref{eqn:low_rank} significantly simplifies the computation of the encoder and decoder functions. This simplification allows for the utilization of only pointwise multiplication and fast Fourier transform operations. As a result, the computation complexity of both the encoder and decoder becomes almost linear, which greatly enhances the efficiency of the pd-IAEnet.

	\subsection{Miscellaneous}
	To enhance the performance of the IAE-net framework, several key features were incorporated, as detailed in \cite{ong2022integral}. Notably, the adoption of multi-channel IAE blocks played a crucial role in extracting additional features from the input data. To tackle common optimization challenges like the vanishing gradient problem, dense skip connections were introduced, connecting individual IAE-blocks. Additionally, a data augmentation process, extensively discussed in \cite{ong2022integral}, was integrated to fully exploit the discretization invariance property, ultimately bolstering the network's generalization capabilities. It is our belief that these enhancements can also prove beneficial for the pd-IAEnet. Consequently, we offer a concise overview of these concepts for the sake of comprehensiveness.

	\subsubsection{Multi-channel pd-IAE blocks}
	A multi-channel pd-IAE block has the following computational flow
	\begin{equation}\label{eqn:multi-channel}
		\begin{tikzcd}
			& a_1 \ar[r,"\text{encoder}"] & v_1 \ar[r,"\text{FNN}"] & u_1 \ar[r,"\text{decoder}"] & b_1  \ar[r,"F_1^{-1}"] &\Tilde{b}_1 \ar[dr]\\
			a \ar[ur,"F_1"] \ar[dr,"F_2"']   &  & & & & & \begin{bmatrix}\Tilde{b}_1 \\ \Tilde{b}_2 \end{bmatrix} \ar[r,"\text{FNN}"]&b\\
			& a_2 \ar[r,"\text{encoder}"] & v_2 \ar[r,"\text{FNN}"] & u_2 \ar[r,"\text{decoder}"] & b_2 \ar[r,"F_2^{-1}"] & \Tilde{b}_2 
			\ar[ur]
		\end{tikzcd}
	\end{equation}
	The input function is duplicated into two separate channels: $a_1 = F_1(a)$ via an identity mapping, denoted as $F_1$. Concurrently, it is transformed into another channel, $a_2 = F_2(a)$, where $F_2$ could be a Fourier or Wavelet transform. 
	Both channels, $a_1$ and $a_2$, are processed individually through a standard pd-IAE block, yielding two distinct outputs, $b_1$ and $b_2$. The output $b_1$ is transformed back into $\Tilde{b}_1 = F_1^{-1}(b_1)$ via an identity mapping, represented as $F_1^{-1}$. Simultaneously, the output $b_2$ is transformed into $\Tilde{b}_2 = F_2^{-1}(b_2)$ using an inverse Fourier or Wavelet transform, denoted as $F_2^{-1}$.  Both transformed outputs, $\Tilde{b}_1$ and $\Tilde{b}_2$, are concatenated into a single extended vector. This concatenated vector is then further processed through a FNN to produce the final output $b$.
	
	It's important to note that the inclusion of the additional Fourier (or Wavelet) channel proves particularly beneficial when working with data featuring oscillatory or sparse characteristics. Essentially, the Fourier channel can be viewed as the application of a pseudo-differential operator (as defined in Equation \eqref{eqn:def_pdo}), followed by an FNN, and another pseudo-differential operator. This two-channel pd-IAE block structure can be readily extended to accommodate scenarios with multiple channels, and it effectively maintains the discretization invariant property.
	
	\subsubsection{Dense skip connections}
	To enhance the input features and mitigate gradient vanishing issues commonly associated with Feedforward Neural Networks (FNNs), skip connections were recursively introduced in each pd-IAE block, a technique inspired by \cite{he2016deep}. These skip connections contribute to greater training stability.
	
	The recursive skip connections augment the basic pd-IAE-net structure \eqref{eqn:basic_IAE} as follows: each intermediate function $a_i$ undergoes an affine transformation $\Acal_i$ and is then added to all subsequent functions $a_j$, where $j>i$. In other words, any intermediate function $a_i$ can be expressed as the sum of the output $\Ical_i(a_{i-1})$ from the preceding IAE-block and a series of affine transformations applied to all preceding functions:
	\begin{equation*}
		a_i = \Ical_i(a_{i-1}) + \sum_{j=0}^{i-1} \Acal_j(a_j)\,, \quad i = 1,\ldots,L \,.
	\end{equation*}
	
	In comparison to the basic pd-IAEnet structure \eqref{eqn:basic_IAE}, the final pd-IAE network takes on the following configuration:
	\begin{equation}\label{eqn:skip_IAE}
		\begin{tikzcd}
			\Bar{f} \ar[r,"F"]  &a_0 \ar[r,"\Ical_1"] \ar[r,bend right, red] \ar[rr, bend right, red] \ar[rrrr, bend right, red]&a_1 \ar[r,"\Ical_2"] \ar[r, bend right, red] \ar[rrr,bend right, red] &a_2 \ar[r,"\Ical_3"] \ar[rr,bend right, red]&\cdots \ar[r,"\Ical_L"] &a_L \ar[r,"G"] &\Bar{g}
		\end{tikzcd}
	\end{equation}
	The red arrows in the diagram signify the skip connections established between different layers through affine transformations. Notably, each pd-IAE block $\Ical_i$ may incorporate multiple channels with dense skip connections, thereby enhancing the network's expressive power.
	
	\subsubsection{Data augmentation training}\label{subsec:da_training}
	
	Neural networks that lack discretization invariance often suffer from overfitting and poor generalization when tested with different discretizations. In contrast, the pd-IAE net benefits from the discretization invariant property, allowing it to be trained and tested with data of varying discretizations while maintaining consistent accuracy. To fully leverage this advantage, the training dataset, initially of a fixed discretization format, is augmented with data pairs representing different resolutions. This augmentation can be easily accomplished using standard interpolation algorithms.
	
	Let's assume that the training data is discretized on grid points:  $S_\Xcal^0 = \{x_1^0,x_2^0,\ldots,x_{s_0}^0\}$ for the input and  $S_\Ycal^0 = \{y_1^0,y_2^0,\ldots,y_{s_0}^0\}$  for the output. We denote $S_\Xcal^i$ and $S_\Ycal^i$, for $i = 1,\ldots,T$, as sequences of different grid points.We can define interpolation operators from $S_0$ to $S_i$ as $\Ical^i_\Xcal: f|_{S_\Xcal^0} \mapsto f|_{S_\Xcal^i} \,, i=1,\ldots,T$ for any input function $f$ and analogously $\Ical^i_\Ycal: g|_{S_\Ycal^0} \mapsto g|_{S_\Ycal^i} \,, i=1,\ldots,T$  for any output function $g$. 
	
	The data augmentation process generates new training data with various discretizations by interpolating the original training data's discretization to a random discretization from the sequences $S_\Xcal^i$ and $S_\Ycal^i$, respectively.

	The loss function of training a pd-IAE net $\Phi_\text{IAE}(\cdot;\theta_{\Phi_\text{IAE}})$ can be defined as the following: 
	\begin{equation*}
		\min_{\theta_{\Phi_\text{IAE}}} \Ebb_{(\Bar{f},\Bar{g})\sim \pi_{\text{data}}} \Ebb_{(I_\Xcal,I_\Ycal) \sim \pi_\text{int} } \left[ L\left(\Phi_\text{IAE}(\Bar{f};\theta_{\Phi_\text{IAE}}),\Bar{g}\right) 
		+ \lambda L\left(\Phi_\text{IAE}(I_\Xcal(\Bar{f});\theta_{\Phi_\text{IAE}}),I_\Ycal(\Bar{g}\right)
		\right]
	\end{equation*}
	where $\pi_\text{data}$ is an independent random measure over the training data set, $\pi_\text{int}$ is an independent random measure over interpolator functions sets $\{\Ical_\Xcal^i\}_{i=1}^T\times\{ \Ical_\Ycal^i\}_{i=1}^T$, $\lambda$ is a hyperparameter that adjusts the balance between augmented data and original data, and $L(\cdot,\cdot)$ is the loss function. In practice, both expectations in the loss are implemented with finitely many empirical data sets and interpolator functions.

	\subsection{Complexity Analysis}
	The proposed pd-IAE net offers a computational advantage over the original IAE-net. Instead of utilizing a dense integration kernel matrix parametrized by neural networks in the IAE-net (cf. \eqref{eqn:iae-en} and \eqref{eqn:iae-de}), the pd-IAE net parametrizes the integral kernel through a low-rank factorization. Consequently, the output sizes are much smaller, resulting in reduced neural network sizes and faster evaluation.

	Assuming that the data matrices  $\Psf,\Qsf,\Tilde{\Psf}$ and $\Tilde{\Qsf}$ are given, the computational complexity of both Algorithm \ref{alg:encoder} and Algorithm \ref{alg:decoder} is $\Ocal(Ks + K s\log s + Km) = \Ocal(Ks\log s)$. It is worth noting that the intermediate  FNN  in the pd-IAE block has an input and output size of $m$. If we assume the FNN has an order $\mathcal{O}(1)$ number of hidden layers with $m$ neurons, then the cost of FNN evaluation is on the order of $\Ocal(m^2)$. Therefore, the total complexity of evaluating a pd-IAE block is $\Ocal(Ks\log s + m^2)$, where $m \leq s$ can be chosen at will. In comparison, the complexity of the original IAE-block in \eqref{eqn:basic_IAE} is at least $\Ocal(ms + m^2)$, even if linear integral transform encoders are used. The pd-IAE block has a computational complexity advantage when $K \ll m \ll s$. 
	
	In the numerical section, we will demonstrate that the low-rank factorization structure of the pd-IAE net not only reduces the number of trainable parameters, leading to faster evaluation, but also results in smaller generalization errors compared to the IAE-net. This is possibly because the smooth assumption in \eqref{eqn:low_rank} for the integral kernel helps mitigate overfitting to the training dataset. 
	Additionally, we will demonstrate that this flexible structure enables the pd-IAE net to have significantly fewer trainable parameters than other benchmark neural networks such as FNO and DeepONet.
	\section{Numerical Experiments}\label{sec:numerics}
	We have conducted a performance analysis of the pd-IAE net alongside several baseline models on benchmark inverse problems, including electricity impedance tomography (EIT), inverse scattering, and optical tomography (OT). Below, we provide brief descriptions of the baseline models:
	\begin{itemize}
		\item \textbf{IAE-net} \cite{ong2022integral}: IAE-net is an integral transform neural network upon which pd-IAE net is built. The primary modification in pd-IAE net is the inclusion of a low-rank factorization (as defined in Equation \eqref{eqn:low_rank}), which results in a smaller neural network size and faster evaluation.
		\item \textbf{Fourier Neural Operator (FNO)} \cite{FNO}: FNO is a novel integral neural network that incorporates Fourier transforms, making it adaptable for training and testing on data with various discretizations. The FNO network consists of a chain of FNO blocks, comprising a Fourier transform, a linear neural network parametrized transformation, and an inverse Fourier transform.
		\item \textbf{DeepONet} \cite{DeepONet}: DeepONet is a versatile operator learning approach, featuring trunk and branch networks. Initially, DeepONet exhibits semi-discretization invariance, with the branch network relying on the input data's size. In our implementation, we have added an interpolation operator to DeepONet, ensuring discretization invariance.
		\item \textbf{ResNet} \cite{he2016deep}: ResNet is a modern neural network that incorporates identity shortcut connections. In our implementation, we use convolutional layers, which naturally exhibit discretization invariance. However, we omit layers that lack discretization invariance, such as a fully connected final layer commonly found in the ResNet architecture. This approach allows ResNet to avoid imposing specific discretization on the input data.
	\end{itemize}
	
	\subsection{Training details}
	Unless otherwise specified, we adhered to consistent training procedures across all benchmark problems. Our data generation process commenced with the creation of 11,000 data pairs, involving target parameters $a(x)$ and corresponding measurements $\Scal_a^n$ (as defined in Equation \eqref{eqn:finitedata}) without introducing any noise. This dataset was then randomly partitioned, with 10,000 data pairs allocated to the training dataset and 1,000 data pairs designated for the test dataset. For simplicity, we assume a uniform computational domain for the target parameter $a(x)$ in all problems, encompassing the unit square $[0,1]^2$, discretized on a uniform mesh denoted as $\{S_\Ycal^0$.
	
	Prior to initiating the training process, we conducted data preprocessing through min-max normalization, which rescaled all input and output data to fit within the range  $[0, 1]$. 
	
	All model training took place on a 48 GB Nvidia Quadro RTX 8000. Due to the fine discretization employed, a uniform batch size of 5 was maintained across all experiments to prevent GPU memory overflow.
	
	For pd-IAE net and IAEnet, we harnessed data augmentation training, as introduced in Section \ref{subsec:da_training}, to maximize efficiency. Conversely, other models did not employ data augmentation training. This data augmentation process involved random discretization selection, followed by the application of bicubic interpolation to adjust both the input and output data to match the chosen discretization.
	
	Evaluation of relative errors across all models followed a consistent methodology, occurring at the conclusion of each epoch. For each of the 1,000 test data points, we interpolated each data pair across all predetermined discretizations. Subsequently, the relative error was computed for each discretization, enabling the determination of an overall average relative error.
	
	The optimization process employed mean squared error as the loss function. All models were trained using the Adam algorithm, without regularization, with a fixed learning rate of 0.001. In the event that the testing average relative error failed to decrease within 40 epochs, the learning rate was halved. Training was halted if there was no improvement in the testing average relative error over 100 epochs.

	For hyperparameter selection, we relied on established configurations from prior research or the corresponding author's code. Hyperparameters and main architecture choices for all experiments were consistent and are detailed as follows:
	\begin{itemize}
		\item \textbf{pd-IAE net:} This network utilized 4 multi-channel pd-IAE blocks in our experiments. Each multi-channel pd-IAE block included an identity channel and a Fourier channel, as shown in Equation \eqref{eqn:multi-channel}. Encoders and decoders were configured with a uniform hyperparameter $m=12$, determining the truncation size of high-frequency components (\ref{alg:encoder}). Additionally, the hyperparameter $K$, representing the rank of the approximation for the kernel function (\ref{eqn:low_rank_encoder2}), was set to 3.
		
		\item \textbf{IAE-net:} Similar to the architecture of pd-IAE net, this network comprised 4 multi-channel IAE blocks, each containing one identity channel and one Fourier transform channel. IAE-net featured a hyperparameter that determined the size of the encoded data after the integral transform within an IAE block, which was set to 12. We obtained the code for this model from the authors of \cite{ong2022integral}. 
		
		\item \textbf{FNO:} This model was configured sequentially, comprising 4 FNO blocks. It featured a hyperparameter defining the number of frequency components to be processed in the Fourier domain, equivalent to the truncation size of high-frequency components in pd-IAE net. We set this value to 12 to match pd-IAE net. The code was implemented based on the one provided in \cite{FNO}.
		
		\item \textbf{DeepONet:} This model included two main components: the trunk net and the branch net. Both the trunk and branch nets were designed as fully connected neural networks with 5 hidden layers, respectively. The trunk net maintained constant hidden layer widths of 256, while the branch net's width progressively decreased from 2048 to 256. These configurations were based on observations from DeepONet's literature and code library available in \cite{DeepONet}.
		
		\item \textbf{ResNet:} This model consisted of 8 residual blocks. Each block incorporated convolutional layers with decreasing kernel sizes, accompanied by a residual skip connection.
		
	\end{itemize}
	
	\begin{table}[htbp]
		\centering
		\begin{tabular}{lrrr}
			\hline
			Model     & Parameters &  Disk space \\
			\hline
			pd-IAEnet &  5,618,785 & 64.8 MB  \\
			IAEnet    & 6,040,657 & 69.4 MB  \\
			FNO       & 9,462,849 & 108 MB   \\
			DeepONet  &  81,597,697 & 933 MB \\
			ResNet    &  5,534,593 & 63.5 MB  \\
			\hline
		\end{tabular}
		\vspace{0.2 cm}
		\caption{This table presents the evaluation time, the number of trainable parameters, and the corresponding disk space requirements for each model. Throughout our experiments, both the number of parameters and disk space usage remained relatively consistent. It is worth noting that pd-IAE net utilizes fewer trainable parameters and requires less disk space compared to IAEnet, FNO, and DeepONet.}
		\label{tab:time_params_disk}
	\end{table}
	\Ke{change resnet to 10 residual blocks}
	In Table \ref{tab:time_params_disk}, we provide an overview of the number of trainable parameters for all models. It is evident that pd-IAEnet outperforms other models in terms of parameter efficiency.
	
	To be more specific, pd-IAE net boasts the fewest parameters among all models. It is approximately 60\% the size of FNO and significantly more compact than DeepONet. It's worth noting that DeepONet's higher parameter count can be attributed to its wider hidden layers. Notably, ResNet features the smallest number of parameters.
	
	We employ the average relative error to assess a model's accuracy and its discretization invariance. This measure is defined as follows:
	\begin{equation}\label{eqn:relaerr}
		\text{relative error} = \frac{\| a - \hat{a}\|_2}{\|a\|_2} \,,
	\end{equation}
	where $\hat{a}$ represents the model's prediction, and $a$ corresponds to the ground truth medium. We calculate the average relative error across all discretizations based on the 1,000 samples in the test dataset. Subsequently, we create plots of the average relative error over each epoch of the training process.

	\subsection{Inverse scattering}
	In the context of the inverse scattering problem presented in Equation \eqref{eqn:helmholtz}, we examined two distinct distributions of scatters, denoted as $\eta$. One distribution featured multiple point scatters, while the other involved Shepp-Logan phantom scatters.

	\subsubsection*{Point scatter distribution}
	For the point scatter case, we generate data using the finite difference method and Neumann series expansion, as shown in \cite{khoo2018switchnet}.The point scatter medium consisted of isolated points with small supports, and we evaluated it using single-frequency measurements. These point scatters were defined within the domain $\Omega = [-0.5, 0.5]^2$ and discretized on an 81x81 grid. The frequency was set to $\omega = 18\pi$, equivalent to a source frequency of 9 Hz. The $\eta(x)$  distribution represented by point scatters comprised four Gaussian mixtures that resembled isolated points. Both the sources and receivers were aligned on a straight line directly above the $\Omega$ domain, and each was discretized into 81 points. Consequently, both the scatter and measurement dimensions were 81x81.

	\begin{figure}[htbp]
		\centering
		\includegraphics[width=0.9\textwidth]{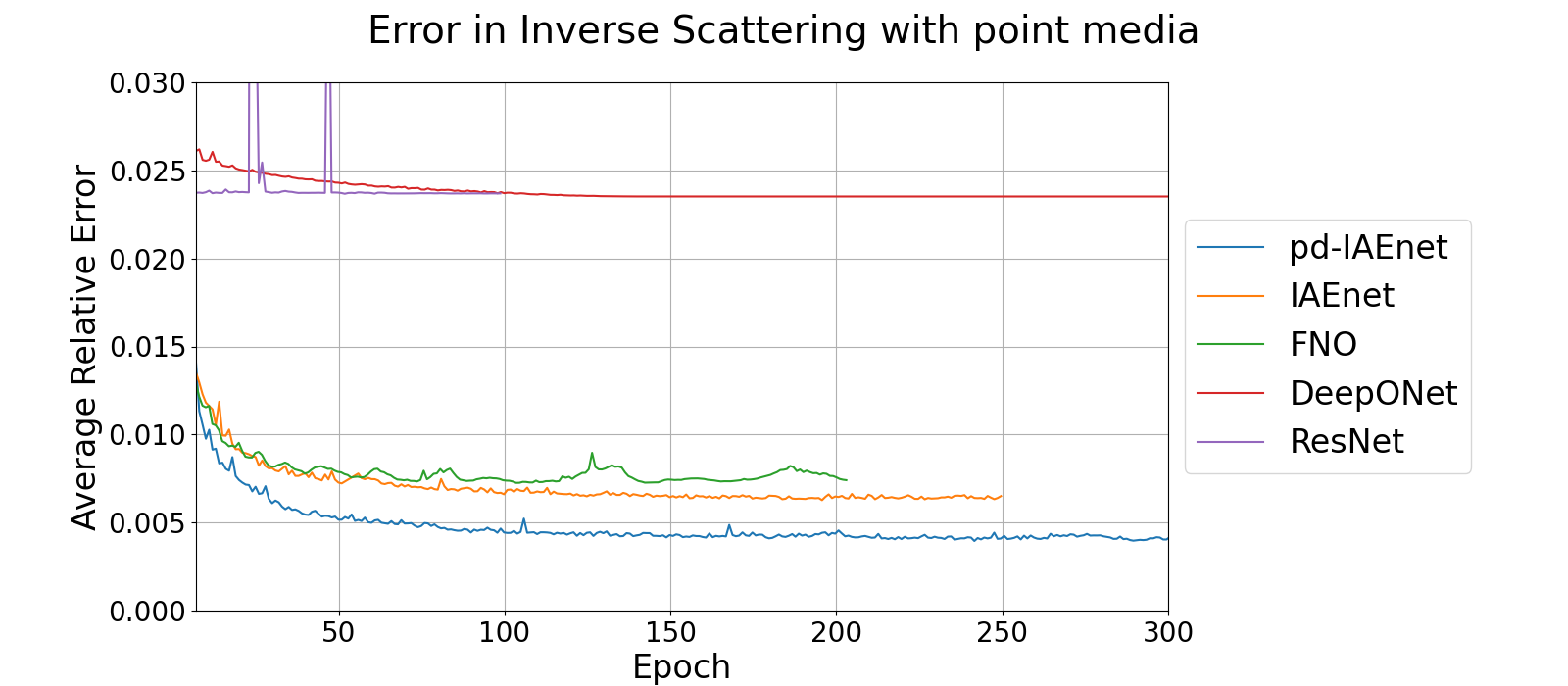}
		\caption{The plot illustrates the average relative error versus the number of epochs for the inverse scattering problem with point media. The y-axis represents the average relative errors. Notably, pd-IAE net consistently achieves the lowest average relative error among all the models.}
		\label{fig:error_point_scattering}
	\end{figure}
	
	In Figure \ref{fig:error_point_scattering}, we calculated the average relative error by averaging errors across various selected discretizations, as defined in \eqref{eqn:relaerr}. We then plotted this error against the training epochs for the inverse scattering problem with point media. The average error was calculated from the following discretizations: $27 \times 27$, $41 \times 41$, $81 \times 81$, $161 \times 161$, and $241 \times 241$.
	It's worth noting that DeepONet and ResNet struggled to learn the inverse operator, while pd-IAE net consistently achieved the lowest error among all baseline models, demonstrating its superior overall accuracy across different discretizations.
	
	\begin{figure}[htbp]
		\centering
		\begin{subfigure}{0.33\textwidth}
			\includegraphics[width=\linewidth]{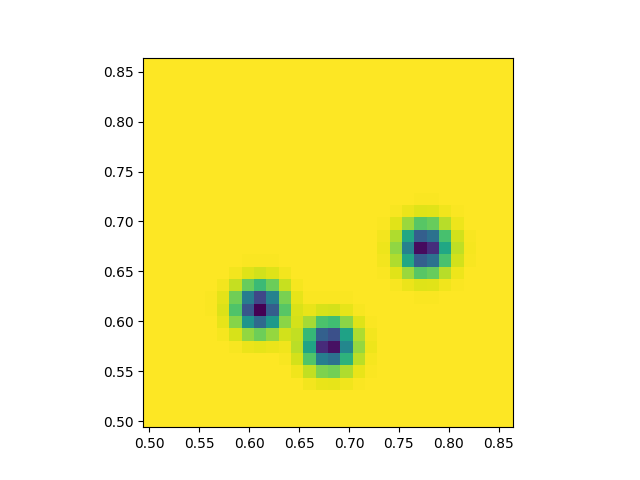}
			\caption{Groundtruth}
		\end{subfigure}%
		\hfill
		\begin{subfigure}{0.33\textwidth}
			\includegraphics[width=\linewidth]{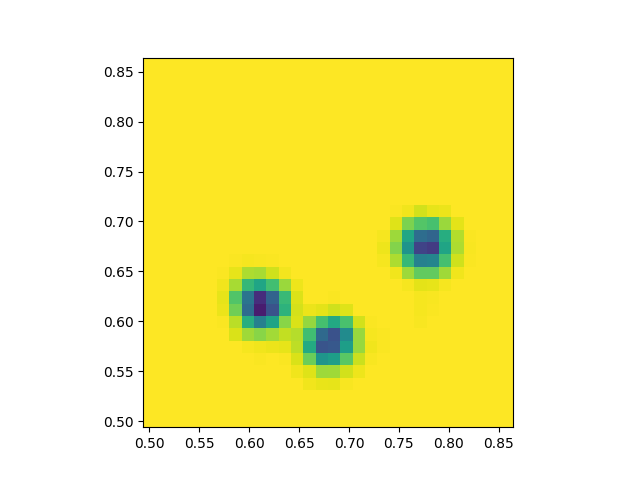}
			\caption{pd-IAE net}
		\end{subfigure}%
		\hfill
		\begin{subfigure}{0.33\textwidth}
			\includegraphics[width=\linewidth]{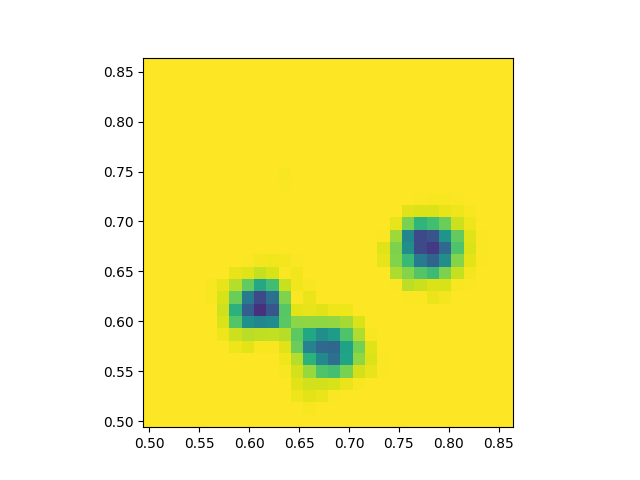}
			\caption{IAEnet}
		\end{subfigure}
		
		
		\begin{subfigure}{0.33\textwidth}
			\includegraphics[width=\linewidth]{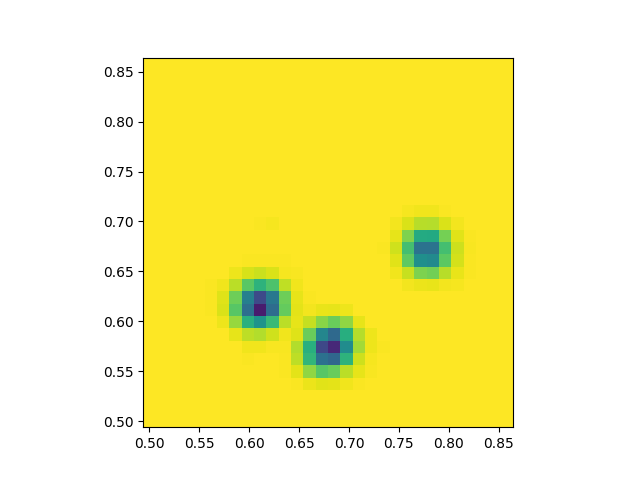}
			\caption{FNO}
		\end{subfigure}%
		\hfill
		\begin{subfigure}{0.33\textwidth}
			\includegraphics[width=\linewidth]{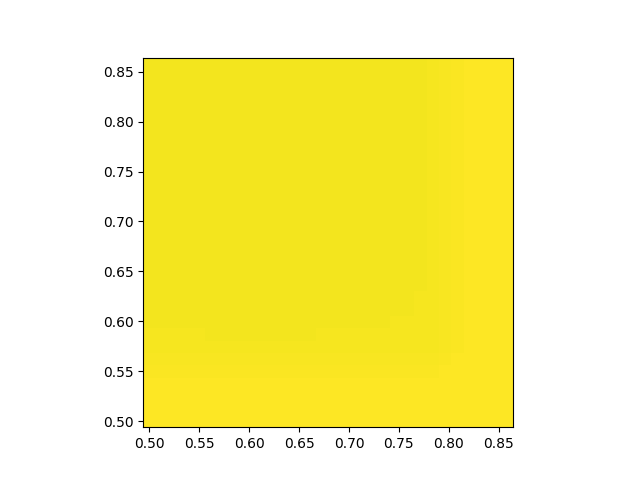}
			\caption{DeepONet}
		\end{subfigure}%
		\hfill
		\begin{subfigure}{0.33\textwidth}
			\includegraphics[width=\linewidth]{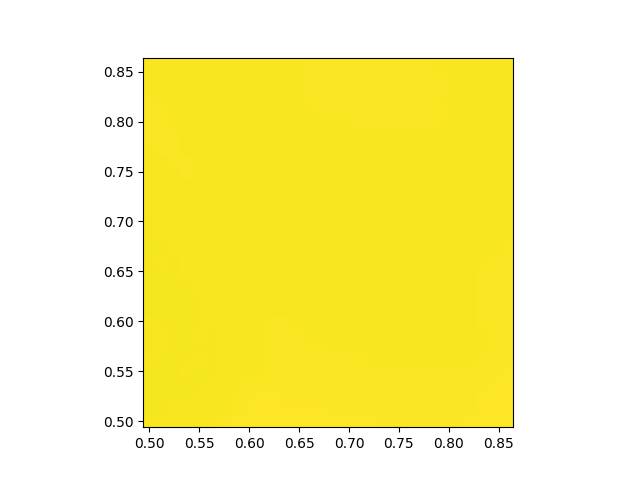}
			\caption{ResNet}
		\end{subfigure}
		
		\caption{In the figures depicting the reconstructed point scatterer at a discretization of $81\times81$, we showcase predictions from all models based on randomly selected test data. Since the point media is uniform across most of the domain, we provide a zoomed-in view of a specific region to emphasize the differences.
			pd-IAEnet accurately reconstructs the points with uniform magnitudes for the three depicted points. In contrast, IAEnet's reconstruction shows the bottom two points merging slightly, despite their separation in the ground truth. FNO, on the other hand, reconstructs the upper-right point with a smaller magnitude compared to the other two points. However, both DeepONet and ResNet struggle to accurately learn the scatterer, with their outputs closely approximating a near-zero solution.}
		\label{fig:inv_scatter_point_reconstruction}
	\end{figure}
	In Figure \ref{fig:inv_scatter_point_reconstruction}, we present sample reconstructed images of the point scatter for each model. These images are discretized on a grid of $81 \times 81$ and rendered to the same scale as the ground truth. Additionally, we zoom in on a specific region to emphasize the differences between pd-IAEnet, IAEnet, and FNO.
	We observed that pd-IAE net, IAEnet, and FNO accurately located the positions of the point scatters. However, IAEnet reconstructed the bottom two points with an overlap not observed in the ground truth. On the other hand, FNO's upper-right point had a noticeably smaller magnitude compared to the other two points. pd-IAE net provided a better-reconstructed image, with the magnitude of each point being more uniform.
	DeepONet and ResNet struggled to reconstruct a meaningful scatter, ultimately producing a near-zero solution. We suspect that the necessity for DeepONet to use interpolation to maintain discretization invariance made the learning more challenging. Additionally, the architecture of ResNet may not have been sophisticated enough to accurately reconstruct small points throughout the domain.

	\subsubsection*{Shepp-Logan scatter distribution}
	\begin{figure}[htbp]
		\centering
		\includegraphics[width=0.9\textwidth]{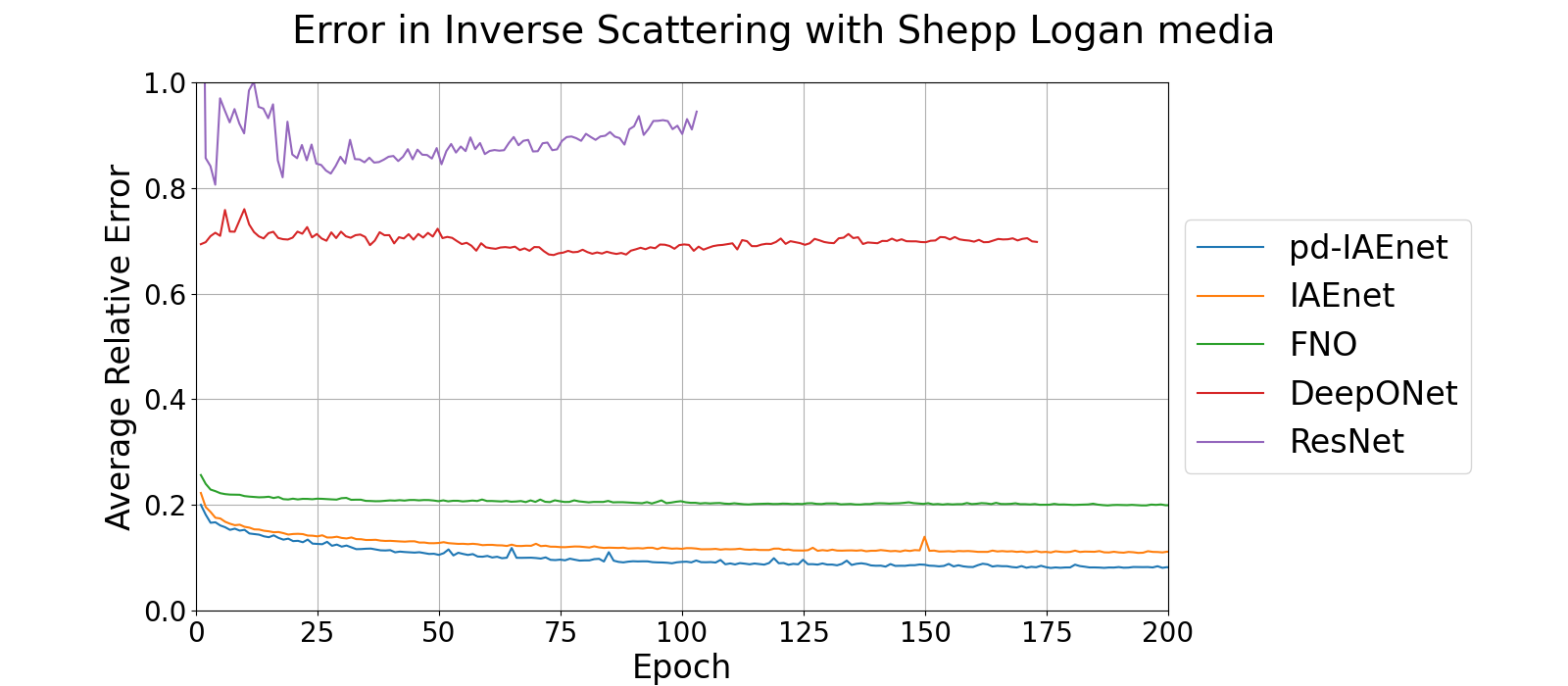}
		\caption{The graph illustrates the average relative error over epochs for the inverse scattering problem with Shepp-Logan media. Notably, among the models, pd-IAEnet consistently achieves the lowest average relative error. 
		}
		\label{fig:error_shepp_logan_scattering}
	\end{figure}
	
	We then proceeded to test the inverse scattering problem using a different scatter distribution known as Shepp-Logan media. This media consists of indicator functions supported on ellipses with varying axis lengths, positions, and rotation angles. We generated the data using the code provided by \cite{zhang2022solving} via the finite difference method. The discretization setup was consistent with that of the point scatter case. Given the increased complexity of Shepp-Logan media compared to point scatters, we explored the domain using multiple source frequencies of $2.5$, $5$, and $10$ Hz, generating multi-frequency measurement data.
	
	In Figure \ref{fig:error_shepp_logan_scattering}, we present the average relative error plotted against the number of epochs for the inverse scattering problem with Shepp-Logan media scatters. The average error is computed from errors obtained at different discretizations: $40 \times 40$, $60 \times 60$, $80 \times 80$, $100 \times 100$, and $240 \times 240$.
	Once more, it is evident that pd-IAE net consistently achieves the lowest relative error when compared to the other models.

	\begin{figure}[htbp]
		\centering
		\begin{subfigure}{0.33\textwidth}
			\includegraphics[width=\linewidth]{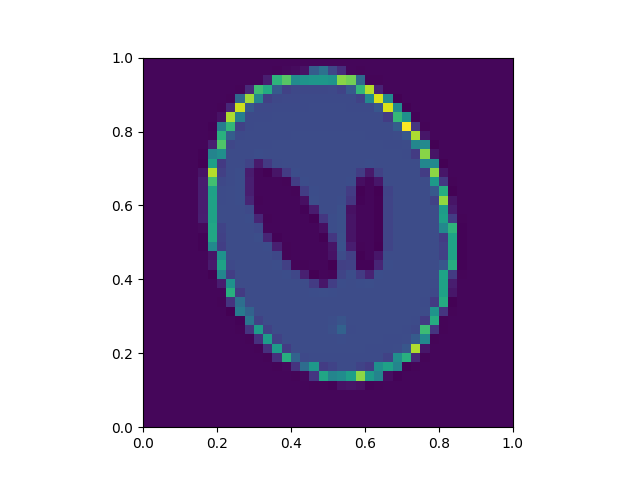}
			\caption{Groundtruth}
		\end{subfigure}%
		\hfill
		\begin{subfigure}{0.33\textwidth}
			\includegraphics[width=\linewidth]{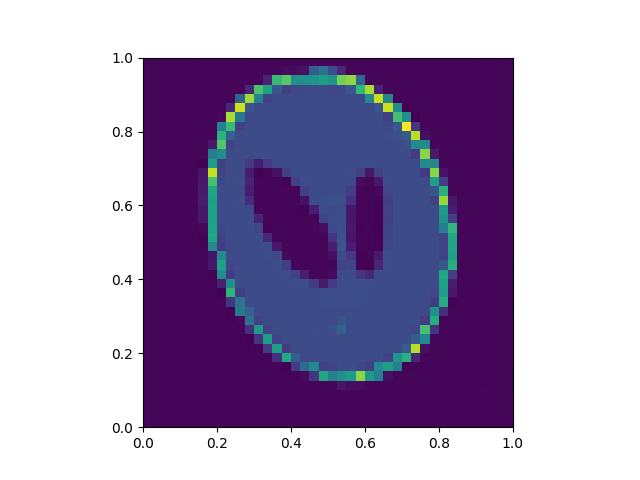}
			\caption{pd-IAE net}
		\end{subfigure}%
		\hfill
		\begin{subfigure}{0.33\textwidth}
			\includegraphics[width=\linewidth]{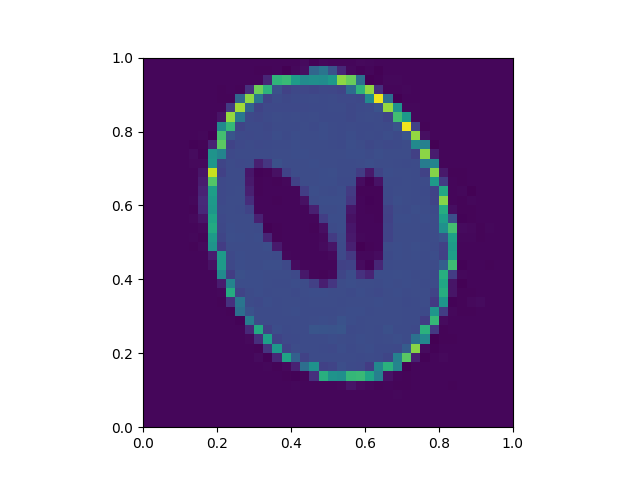}
			\caption{IAEnet}
		\end{subfigure}
		
		
		\begin{subfigure}{0.33\textwidth}
			\includegraphics[width=\linewidth]{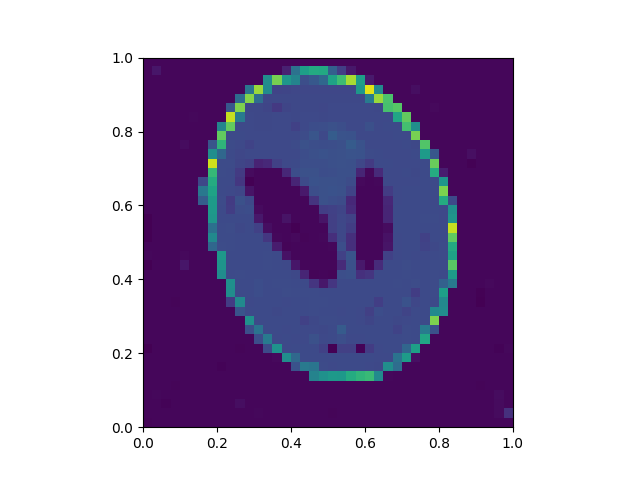}
			\caption{FNO}
		\end{subfigure}%
		\hfill
		\begin{subfigure}{0.33\textwidth}
			\includegraphics[width=\linewidth]{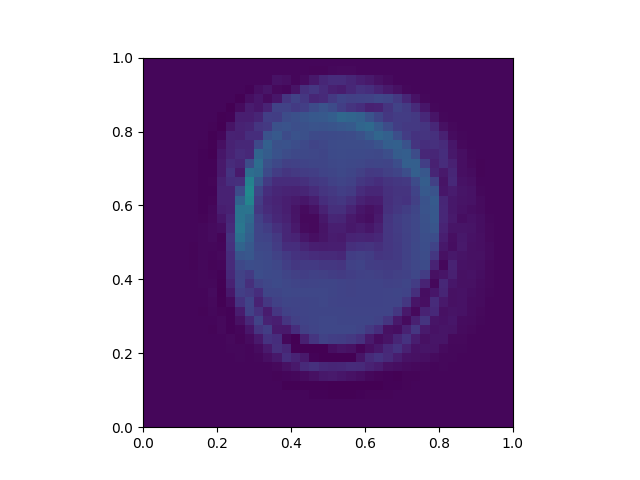}
			\caption{DeepONet}
		\end{subfigure}%
		\hfill
		\begin{subfigure}{0.33\textwidth}
			\includegraphics[width=\linewidth]{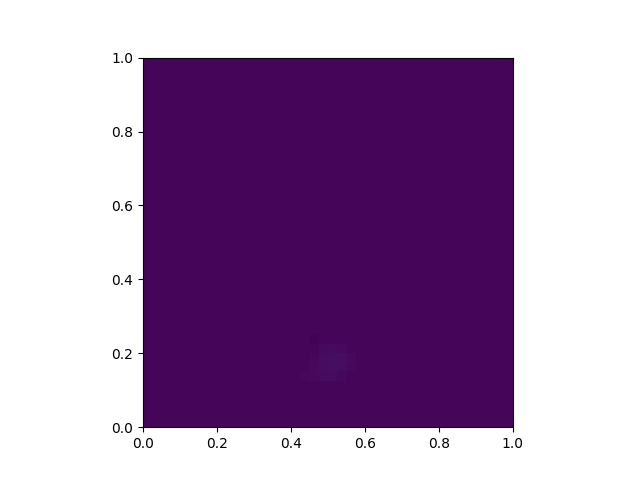}
			\caption{ResNet}
		\end{subfigure}
		\caption{These plots present samples of the reconstructed media at the $40\times40$ discretization level for the inverse scattering problem with Shepp-Logan media. Notably, among the models, pd-IAE net provides the most accurate reconstruction of the ground truth. IAEnet exhibits minor inaccuracies below the two ellipses and along the bottom edge. FNO is marred by small artifacts throughout the image. DeepONet yields a blurry reconstruction with an incorrectly shaped boundary, and ResNet approximates a near-zero solution at the $40\times40$ discretization level.} 
		\label{fig:inv_scatter_shepp_logan_reconstruction}
	\end{figure}
	
	We have also included the reconstructed images for each model using randomly chosen measurement data in Figure \ref{fig:inv_scatter_shepp_logan_reconstruction}. Although the original data is discretized at $80\times80$, we evaluated the models with an input discretized at $40\times40$ to assess their discretization invariance at lower resolutions. We observed notable differences between the reconstruction by pd-IAE net and the ground truth. IAEnet provides a generally good overall reconstruction but introduces a minor inaccuracy below the two internal ellipses and yields a slightly different magnitude along the bottom boundary compared to the ground truth. While FNO captures the overall shape correctly, it exhibits several artifacts, with the most noticeable ones just above the bottom boundary. DeepONet produces a blurry reconstruction of the boundary and the two internal ellipses. In contrast, ResNet struggles at the $40\times40$ discretization level and approximates a near-zero solution.
	
	\subsection{Optical tomography}
	For optical tomography, we employ the finite difference method to solve the radiative transfer equation \eqref{eqn:rte} in a slab plane geometry, as previously discussed in \cite{chen2018stability}. In this configuration, the spatial domain is $[0, 1]$ and the velocity domain is $[-1, 1]$.
	Both the spatial and velocity domains are discretized with $50$ points. The incoming boundary encompasses two spatial endpoints associated with $25$ incoming velocities. Consequently, the measurement matrix is of size $50\times50$. We consider a scattering coefficient, denoted as $\sigma_s(x)$, following a distribution that generates sinusoidal waves with random phase and magnitude.

	\begin{figure}[htbp]
		\centering
		\includegraphics[width=0.9\textwidth]{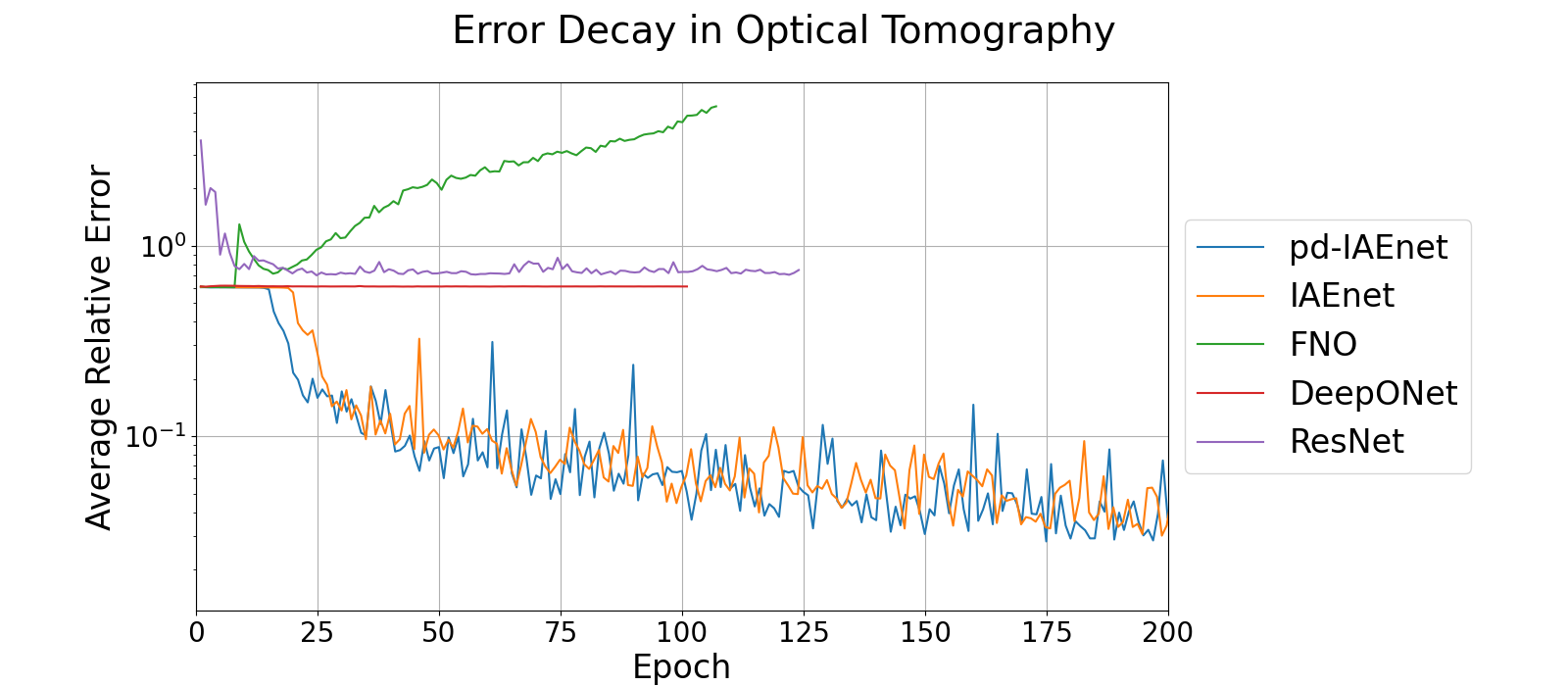}
		\caption{The plot illustrates the average relative error as a function of training epochs for optical tomography. Both pd-IAEnet and IAEnet consistently achieve low errors, each below 0.1. In contrast, DeepONet, ResNet, and FNO encounter challenges in learning the problem in a discretization-invariant manner. DeepONet and ResNet stabilize at a high error level, while FNO exhibits a divergent average relative error.}
		\label{fig:error_RTE}
	\end{figure}
	
	In Figure \ref{fig:error_RTE}, the plot displays the average relative error as a function of training epochs for optical tomography. This average error is computed from the errors at different discretizations: $30 \times 30$, $40 \times 40$, $50 \times 50$, $60 \times 60$, and $70 \times 70$. Since optical tomography is an ill-posed inverse problem with no additional numerical regularization, some methods achieve early convergence with lower accuracy.
	
	Both pd-IAE net and IAEnet perform equally well, with the lowest average relative errors. In contrast, FNO, DeepONet, and ResNet exhibit higher average relative errors. The high error in FNO can be attributed to its decreased accuracy when dealing with discretizations that differ from the original $50\times 50$, as further exemplified in \ref{fig:discretization_RTE}.
	
	\begin{figure}[htbp]
		\centering
		\begin{subfigure}{0.33\textwidth}
			\includegraphics[width=\linewidth]{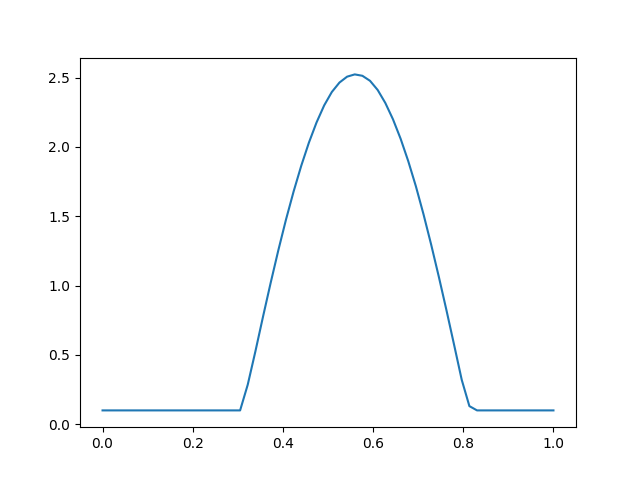}
			\caption{Groundtruth}
		\end{subfigure}%
		\hfill
		\begin{subfigure}{0.33\textwidth}
			\includegraphics[width=\linewidth]{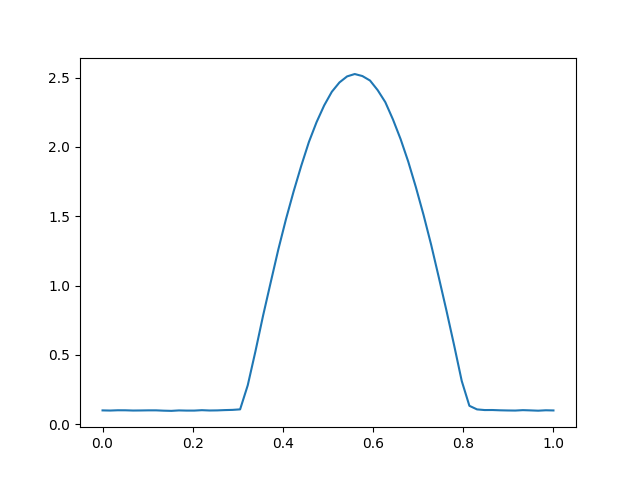}
			\caption{pd-IAE net}
		\end{subfigure}%
		\hfill
		\begin{subfigure}{0.33\textwidth}
			\includegraphics[width=\linewidth]{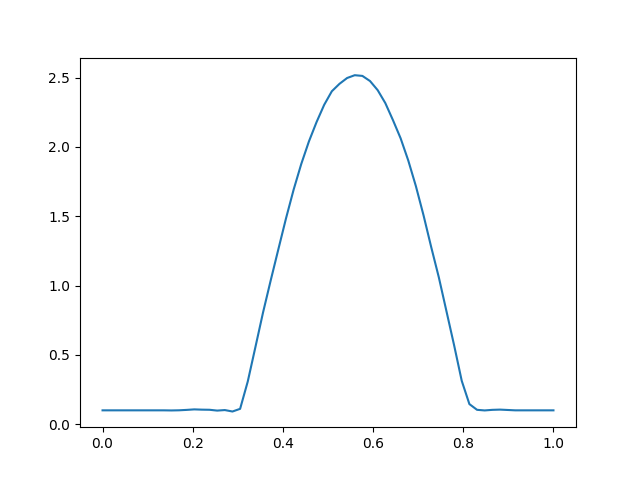}
			\caption{IAEnet}
		\end{subfigure}
		
		
		\begin{subfigure}{0.33\textwidth}
			\includegraphics[width=\linewidth]{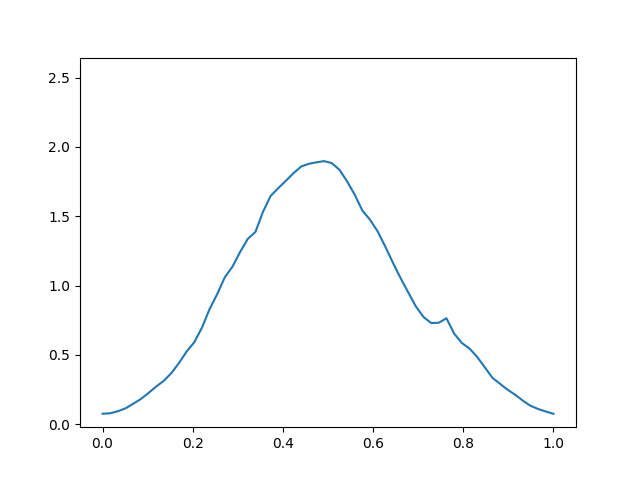}
			\caption{FNO}
		\end{subfigure}%
		\hfill
		\begin{subfigure}{0.33\textwidth}
			\includegraphics[width=\linewidth]{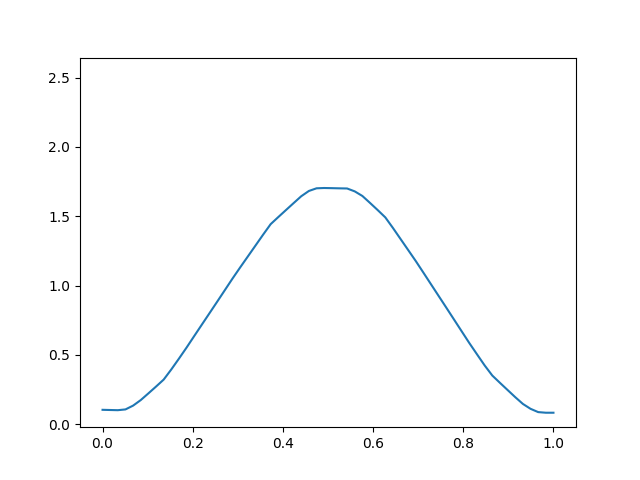}
			\caption{DeepONet}
		\end{subfigure}%
		\hfill
		\begin{subfigure}{0.33\textwidth}
			\includegraphics[width=\linewidth]{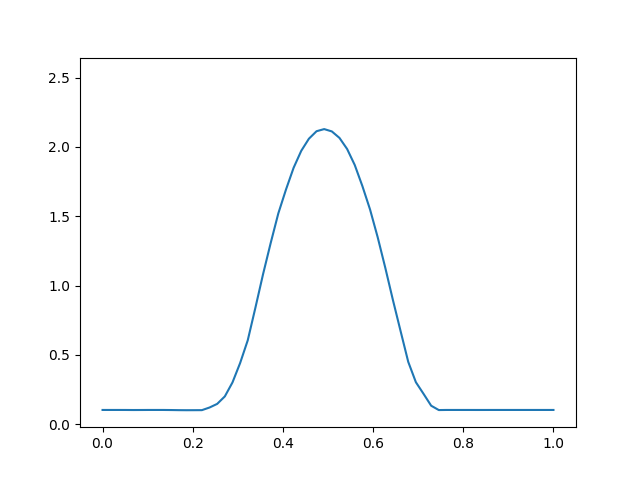}
			\caption{ResNet}
		\end{subfigure}
		
		\caption{The visual representations depict the reconstruction of the sinusoidal scattering coefficient in the context of optical tomography, utilizing a 60x60 discretization grid. Notably, both pd-IAE net and IAEnet yield reconstructions that closely align with the ground truth. However, it is worth noting that IAEnet's reconstruction exhibits minor perturbations, particularly on the left side of the wave's base.
			In contrast, both FNO and DeepONet encounter difficulties in accurately capturing the characteristics of the sine wave, primarily due to their method of averaging all training samples. On the other hand, ResNet excels in identifying the location of the wave but faces challenges in accurately determining its magnitude, revealing a limitation in its reconstruction capabilities.}
		\label{fig:inv_RTE_reconstruction}
	\end{figure}
	
	In Figure \ref{fig:inv_RTE_reconstruction}, we present the reconstructed sinusoidal scattering coefficient for the optical tomography problem with a discretization of $60\times60$. Both Pd-IAEnet and IAEnet deliver reconstructions that closely match the ground truth. Conversely, FNO and DeepONet fall short in accurately capturing the location and wavelength. ResNet, although capable of identifying the location, encounters challenges in determining the correct magnitude of the wave.

	\subsection{Calder\'on problem}
	
	In the context of the Calder\'on problem, we address the elliptic equation \ref{eqn:calderon} through the application of finite element methods. In this study, we consider the unknown medium conductivity, denoted as $a(x)$, which is represented by random Shepp-Logan phantoms discretized on a $64\times64$ mesh covering the unit square $\Omega=[0, 1]^2$. To collect measurements, sources and receivers are strategically positioned along all four sides of the square, resulting in a measurement matrix with dimensions of $252\times 252$.
	\begin{figure}[htbp]
		\centering
		\includegraphics[width=0.9\textwidth]{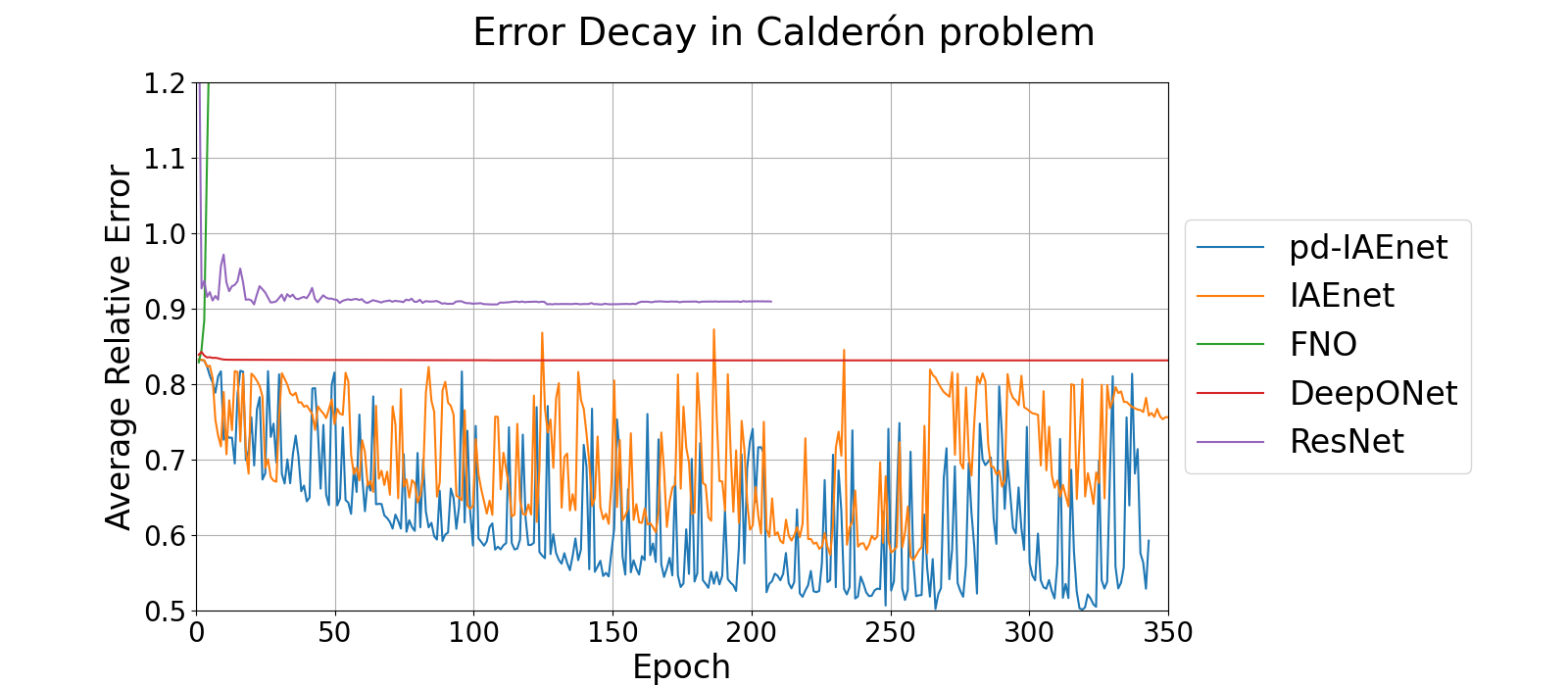}
		\caption{The plot presents the average relative error plotted against the number of epochs for the Calder\'on problem. Notably, FNO exhibits a diverging trend in its average relative error, indicating instability during training. On the other hand, DeepONet and ResNet converge, but with persistently high relative errors. 
			In contrast, both pd-IAE net and IAE-net consistently achieve significantly lower errors, highlighting their robust training. Specifically, IAE-net reaches its smallest error of $56.68\%$  during the training process, while pd-IAE net outperforms with the smallest error of $50.13\%$. }
		\label{fig:error_calderon}
	\end{figure}
	
	In Figure \ref{fig:error_calderon}, we have plotted the average relative error for each model in the context of the Calder\'on problem. The relative error values are averaged across five different discretizations, namely $42\times42$, $63\times63$, $84\times84$, $126\times126$, and $252\times252$. To augment the data, we employed interpolation techniques for both the medium conductivity and measurement data, which were originally provided at discretizations of $64\times64$ and $252\times252$, respectively.
	It's important to note that the relative error metrics for the Calder\'on problem are generally higher compared to other benchmark examples, such as inverse scattering and optical tomography. This is due to the inherently ill-posed nature of the Calder\'on problem. Nevertheless, pd-IAE net stands out by consistently achieving a relative error that is at least $5\%$ smaller than that of all other models. Given the oscillations observed in the testing loss in Figure \ref{fig:error_calderon}, we implemented early stopping during the training of each model to ensure the attainment of the smallest error. This approach allows us to optimize the performance of these models in addressing the challenging Calder\'n problem.
	
\begin{figure}[htbp]
		\centering
		\begin{subfigure}{0.33\textwidth}
			\includegraphics[width=\linewidth]{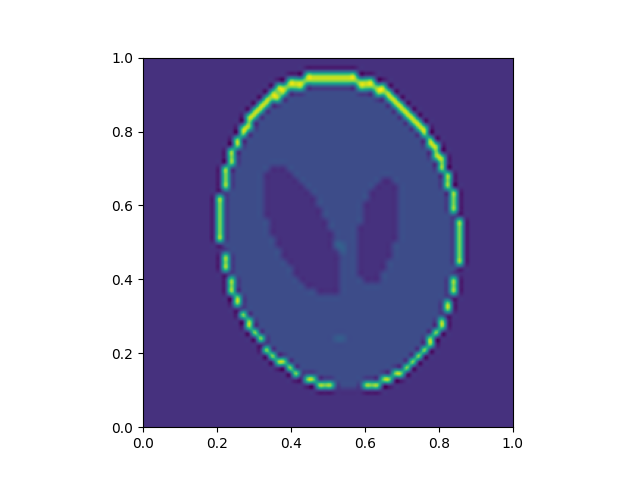}
			\caption{Groundtruth}
		\end{subfigure}%
		\hfill
		\begin{subfigure}{0.33\textwidth}
			\includegraphics[width=\linewidth]{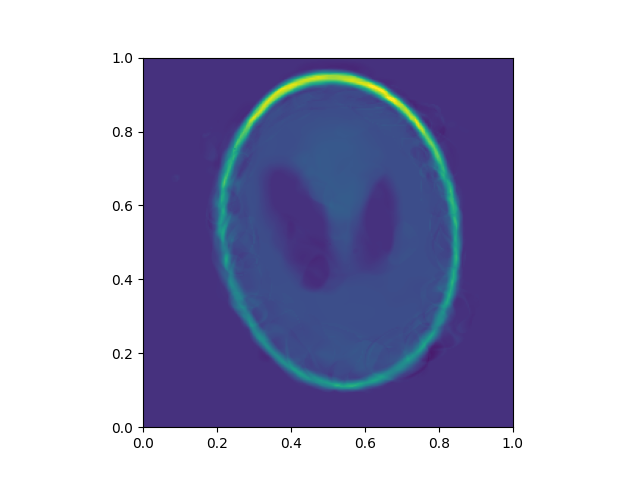}
			\caption{pd-IAE net}
		\end{subfigure}%
		\hfill
		\begin{subfigure}{0.33\textwidth}
			\includegraphics[width=\linewidth]{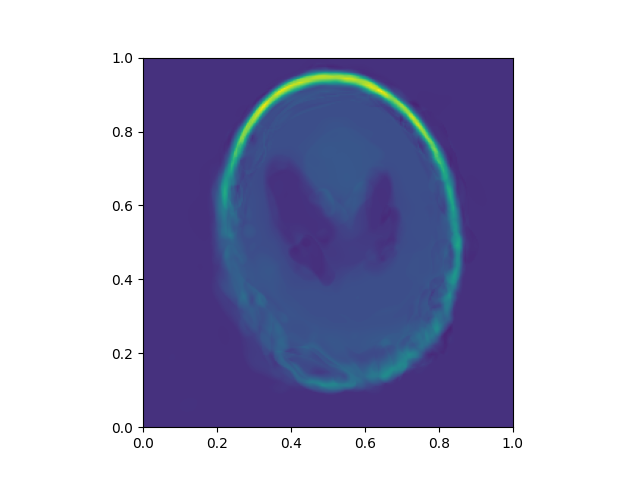}
			\caption{IAEnet}
		\end{subfigure}
		
		
		\begin{subfigure}{0.33\textwidth}
			\includegraphics[width=\linewidth]{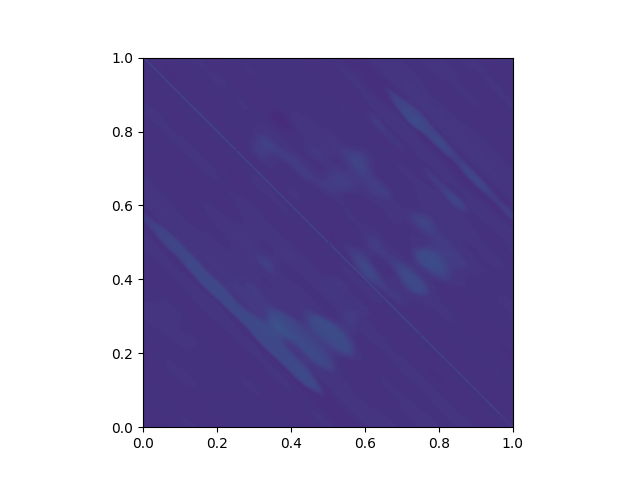}
			\caption{FNO}
		\end{subfigure}%
		\hfill
		\begin{subfigure}{0.33\textwidth}
			\includegraphics[width=\linewidth]{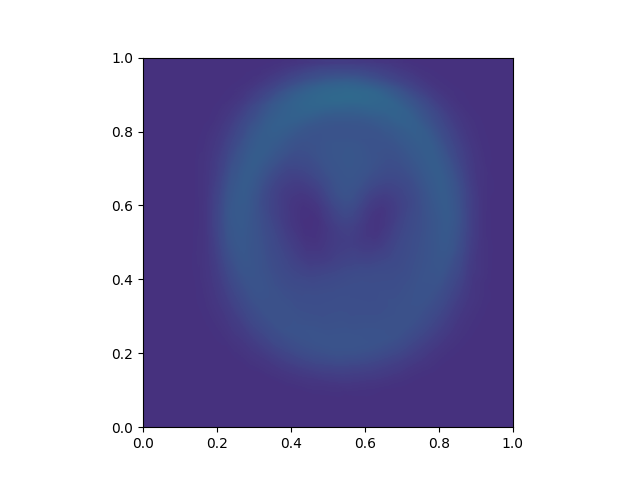}
			\caption{DeepONet}
		\end{subfigure}%
		\hfill
		\begin{subfigure}{0.33\textwidth}
			\includegraphics[width=\linewidth]{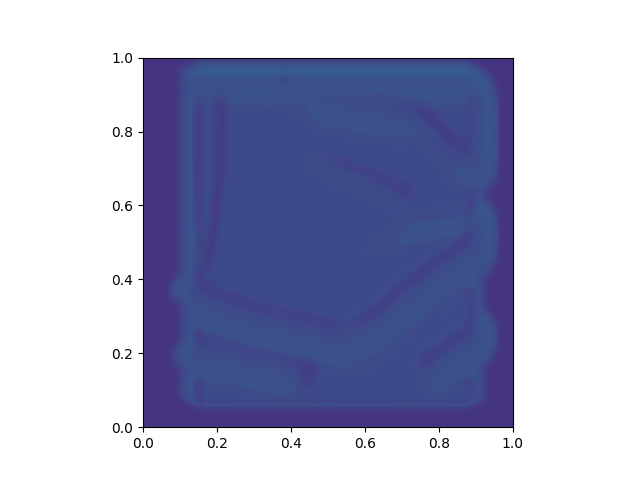}
			\caption{ResNet}
		\end{subfigure}
		
		\caption{The figures display reconstructed images generated by all models based on randomly selected measurement data with a discretization of $252\times252$ for the Calderón problem.
			DeepONet is capable of producing a somewhat blurry image of the Shepp Logan phantom. In contrast, the outputs from FNO and ResNet clearly indicate their inability to accurately reconstruct the Shepp Logan phantom.
			While pd-IAEnet and IAEnet may not capture the intricate details of the phantom's boundary with absolute precision, they consistently produce significantly clearer images when compared to the other models.
		}
		\label{fig:calderon_reconstruction}
	\end{figure}

	In Figure \ref{fig:calderon_reconstruction}, we have selected a random set of measurement data and generated reconstructed images for all models (utilizing early stopping) based on this specific measurement.
	Notably, FNO and ResNet fail to produce reasonable reconstructions, while DeepONet tends to generate an output that is essentially an average of all training media conductivities. In contrast, both pd-IAEnet and IAEnet consistently deliver superior reconstructions when compared to the other models. Of particular significance is the performance of pd-IAEnet, which excels in distinguishing the internal ellipses and provides a more accurate reconstruction of the boundary compared to IAEnet. This observation underscores the superior image reconstruction capabilities of pd-IAEnet in the context of the Calder\'on problem.

	\subsection{Discretization invariance}
	We conducted accuracy tests for all models, assessing their performance across different discretizations. For models that lack full discretization invariance, we utilized naive interpolations to facilitate evaluation. These experiments encompassed the inverse scattering problem involving Shepp-Logan media and optical tomography featuring a sinusoidal scattering coefficient.
	
\begin{figure}[htbp]
		\centering
		\includegraphics[width=0.8\textwidth]{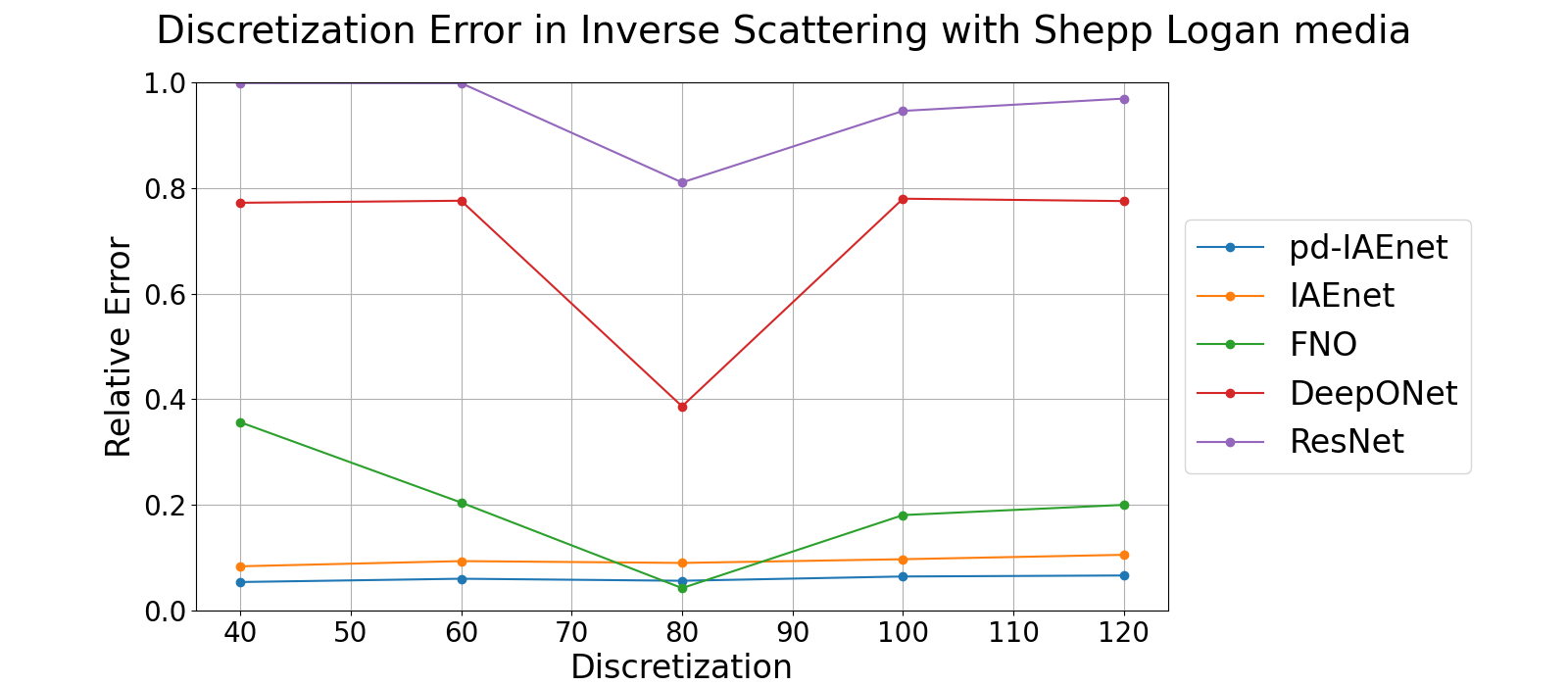}
		\caption{This graph provides insights into the relative error across various discretizations. The x-axis enumerates different discretization sizes: $40\times40$, $60\times60$, $80\times80$, $100\times100$, and $120\times120$. As the discretizations move away from the original setting of $80\times80$, the accuracy of FNO, DeepONet, and ResNet  noticeably declines. In contrast, pd-IAE net and IAEnet consistently maintain high accuracy levels across all discretizations.} 
		\label{fig:discretization_scattering_shepp_logan}
	\end{figure}

	In Figure \ref{fig:discretization_scattering_shepp_logan}, we present the relative reconstruction errors of all models for inverse scattering across a range of discretizations, spanning from $40\times40$ to $120\times120$. Notably, both pd-IAE net and IAEnet demonstrate the remarkable ability to maintain consistent accuracy across all discretizations, highlighting their robustness.
	
	Conversely, FNO stands out for achieving its smallest error on the $80\times80$ mesh, with the error doubling on other discretizations. For DeepONet and ResNet, a similar error bias persists, characterized by significantly larger contrasts in error across the different discretizations.

\begin{figure}[htbp]
		\centering
		\includegraphics[width=0.9\textwidth]{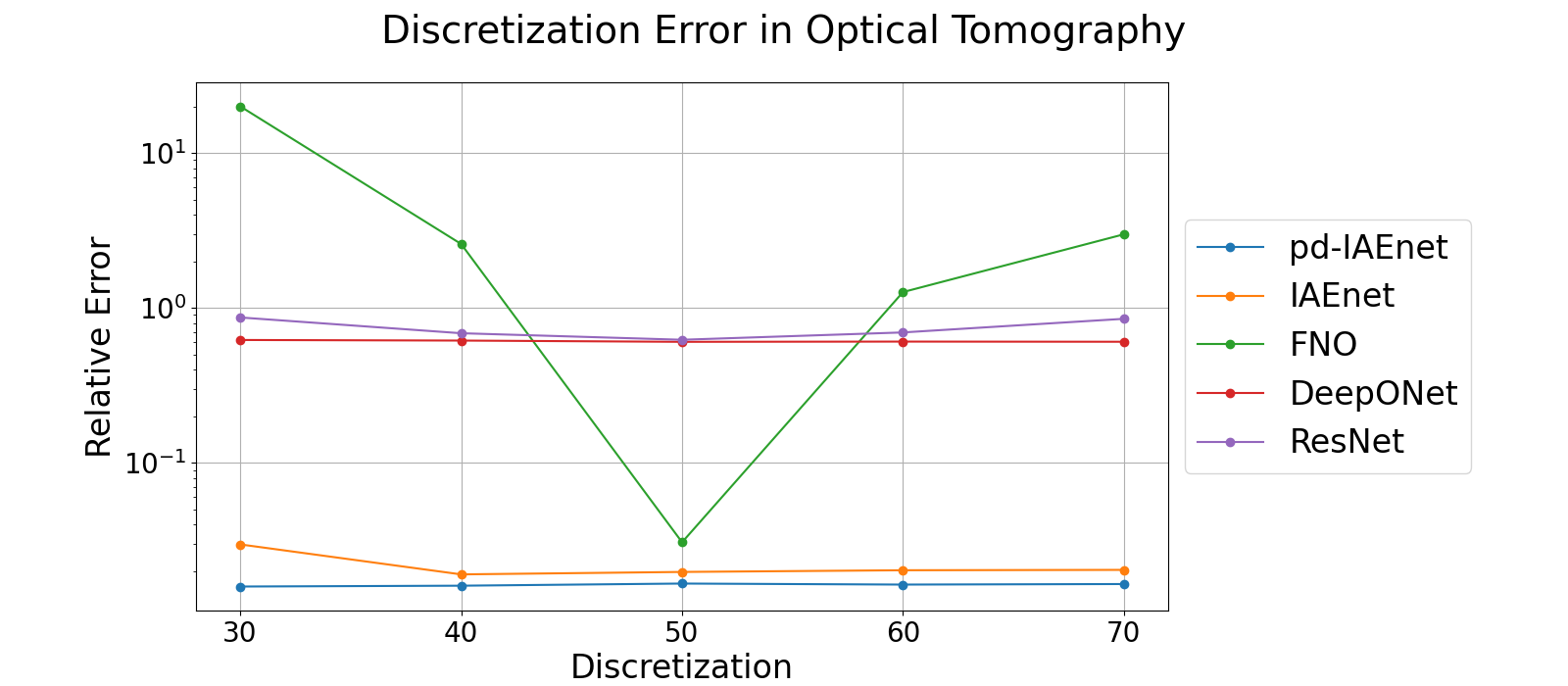}
		\caption{The graph offers a visual depiction of the relative error across a range of discretizations. On the x-axis, you can see the different discretization sizes: $30\times30$, $40\times40$, $50\times50$, $60\times60$, and $70\times70$. pd-IAE net stands out as the model with the best overall accuracy and, notably, it maintains near-uniform accuracy levels across all discretizations. This highlights the model's remarkable consistency and robustness across varying mesh sizes.}
		\label{fig:discretization_RTE}
	\end{figure}

	In Figure \ref{fig:discretization_RTE}, we present the relative error of the reconstructed scattering coefficient in the context of the optical tomography problem. Once more, pd-IAE net and IAE net consistently exhibit uniform accuracy across all discretization meshes.
	FNO, on the other hand, manages to achieve good accuracy on the $50\times50$ mesh but experiences a substantial relative error exceeding $100\%$ on other meshes. In contrast, both DeepONet and ResNet display uniform errors across all discretization meshes, albeit with lower overall accuracy levels.
	
	\subsection{Robustness to noise}
	The ill-posed nature of inverse problems often amplifies noise in the reconstructed parameters. To investigate the sensitivity of reconstruction parameters to noise in the measurement data, we introduced noise to both the training and testing data. Consequently, all models underwent training and testing with noisy data.
	
	Our experiments covered three benchmark examples, including the inverse scattering problem with point scatters and Shepp-Logan media (with $1\%$ additive noise), as well as the inverse  RTE problems (with $0.25\%$ additive noise). However, for the severely ill-posed Calder\'on problem, we opted not to conduct noise tests, as most models struggled to converge in the noise-free case, as evidenced in Figure \ref{fig:error_calderon}.

	In Table \ref{fig:table_no_noise} and Table \ref{fig:table_noise}, we calculated the relative error in both the noise-free and noisy cases. It was observed that, for the inverse scattering problem with point media, all methods exhibited slightly larger reconstruction errors when trained on noisy data. However, for the Shepp-Logan media, the error introduced by noise was more significant.
	
	In the case of the optical tomography problem (denoted as RTE), the reconstruction error for pd-IAEnet and IAE-net changed notably, ranging from approximately $1\%$ to $40\%$, underscoring the ill-posed nature of the problem. Conversely, the error of DeepONet remained relatively stable for the RTE example, primarily because its reconstruction relies on taking the average of all possible training samples, as seen in Figure \ref{fig:inv_RTE_reconstruction}. For FNO and ResNet, which already exhibited larger errors in the noise-free case, the introduction of noise led to changes in error ranging from $10\%$ to $20\%$.

	\begin{table}[htbp]
		\centering
		\begin{tabular}{lrrr}
			\hline
			Model     &   Scattering (point media) &  Scattering (Shepp-Logan media) &  RTE \\
			\hline
			pd-IAEnet &  0.004679 & 0.09102 &  0.01630\\
			IAEnet    &  0.005768 & 0.09139 & 0.0197\\
			FNO       &  0.007085 & 0.1969 & 0.6059\\
			DeepONet  &  0.02351  & 0.6730 & 0.6090\\
			ResNet    &  0.02372 & 0.8105  & 0.6987\\
			\hline
		\end{tabular}
		\vspace{0.2 cm}
		\caption{The table shows the average relative error among models with no noise. 
		}
		\label{fig:table_no_noise}
	\end{table}
	
	\begin{table}[htpb]
		\centering
		\begin{tabular}{lrrr}
			\hline
			Model     &   Scattering (point media) &  Scattering (Shepp-Logan media) &  RTE \\
			\hline
			pd-IAEnet &  0.005187 & 0.1231 & 0.4042 \\
			IAEnet    &  0.007474 & 0.1212 &  0.4522 \\
			FNO       &  0.007357 & 0.2086 & 0.7704 \\
			DeepONet  &  0.02352  & 0.7150 & 0.6126 \\
			ResNet    &  0.02397 & 0.8804  & 0.7693 \\
			\hline
		\end{tabular}
		\vspace{0.2 cm}
		\caption{The table provides an overview of the average relative error across various models for all problems, following the introduction of noise to the input data. In the case of the inverse scattering problem, we introduced $1\%$ additive noise to the measurement data for both the point media and Shepp-Logan media scenarios. For the optical tomography problem, characterized by its ill-posed nature, a smaller $0.25\%$ noise was added to the measurement data.  
		}
		\label{fig:table_noise}
	\end{table}
	
	\section{Conclusions}
	In conclusion, our study presents the "pd-IAE net" as a novel and promising neural network framework for addressing complex inverse problems. By leveraging the principles of the IAE-net framework and drawing inspiration from pseudo-differential operators, we have developed a model that excels in several key aspects.
	Our comprehensive numerical investigation across diverse benchmark inverse problems has consistently demonstrated the superiority of the pd-IAE net. It achieves higher accuracy while maintaining efficiency, with fewer parameters compared to established baseline models. Furthermore, the model's remarkable consistency in accuracy across different discretizations is a significant advantage, ensuring its adaptability to various real-world scenarios.
	The most notable achievement of the pd-IAE net is its ability to produce robust reconstructions, even when confronted with severely ill-posed inverse problems. This resilience sets it apart from other models that struggle to create meaningful parameter reconstructions under similar conditions. For severely ill-posed inverse problems, pd-IAEnet is able to generate more accurate reconstructions while other models may suffer from blurry boundaries, or not be able to construct a reasonable image.
	
	We notice that noises may lead to large reconstruction error for severely ill-posed inverse problems like the optical tomography problem and the Calder\"on problem. To mitigate this issue, prior information of the target parameter is needed.
	In the future work, we will be concerned about the potential for further improvements in reconstruction quality through the incorporation of numerical regularization techniques during training for high ill-posed inverse problems like Calder\"on problem. This avenue represents a promising direction for our future research endeavors, with the ultimate goal of enhancing the practical applicability and impact of the pd-IAE net framework in addressing inverse problems across diverse fields.
	
	\section*{Acknowledgements} 
	C. W. was partially supported by National Science Foundation Award DMS-2136380 and DMS-2206332.  
	
	\bibliographystyle{plain}  
	\bibliography{references}

\end{document}